\documentclass[10pt,a4paper,oneside,british]{amsart}

\usepackage[margin=1em, labelfont=bf,
                singlelinecheck=true, font=sl, format=hang, center]{caption}

\usepackage{amsmath} 
\usepackage{amsthm}
\usepackage{upref}
\usepackage{enumitem}

\setlist{parsep=5pt, topsep=5pt}

\usepackage[latin1]{inputenc}
\usepackage[british]{babel} 
\usepackage[sc]{mathpazo}
\usepackage{bm}

\makeatletter                         
\g@addto@macro\bfseries{\boldmath}    
\makeatother                          

\newcommand{\ttbinom}[2]{\bigl(\substack{#1\\[1.3pt]#2}\bigr)}

\DeclareMathOperator{\peak}{peak}
\DeclareMathOperator{\pattern}{pat}
\DeclareMathOperator{\des}{des}
\DeclareMathOperator{\asc}{asc}

\newtheorem{thm}{Theorem}[section]
\newtheorem{deff}[thm]{Definition}
\let\Definition=\deff                    
\renewcommand{\deff}{\Definition \rm}    
\newtheorem{lem}[thm]{Lemma}
\newtheorem{cor}[thm]{Corollary}
\newtheorem{examp}[thm]{Example}
\let\Example=\examp                  
\renewcommand{\examp}{\Example \rm}  

\newenvironment{bevis}
{\begin{proof}[\rm \bf Proof.]}
{\end{proof}}

\newcommand{\symg}[1]{S_{#1}}
\newcommand{\fall}[2]{#1^{\underline{#2}}}
\newcommand{\Mset}[3]{{\it M}($#1$,$#2$,$#3$)}
\newcommand{\fourset}[1]{C_{\lambda,#1}}
\newcommand{\chivek}{\boldsymbol{\chi}_7}
\newcommand{\patt}[1]{\pattern_{#1}}
\newcommand{\subsetpn}{D_{n}}
\newcommand{\vecsp}{,\;}
\newcommand{\indfn}{v}

\title{Pattern containment in random permutations}

\author{Jonna Gill}

\address{Department of Mathematics, Link\"oping University, SE-581 83 Link\"oping, Sweden}

\email{jonna.gill@liu.se}

\begin{document}
	
	\begin{abstract}
		This paper studies permutation statistics that count occurrences 
			of patterns. Their expected values on a product
			of $t$ permutations chosen randomly from $\Gamma \subseteq S_{n}$, 
			where $\Gamma$ is a union of conjugacy classes, are considered.
			Hultman has described a method for computing such an expected value, 
			denoted $\mathbb{E}_{\Gamma}(s,t)$, of a statistic $s$,
			when $\Gamma$ is a union of conjugacy classes of $S_{n}$. The
			only prerequisite is that the mean of $s$
			over the conjugacy classes is written as a linear combination
			of irreducible characters of $S_{n}$.
			Therefore, the main focus of this article is to express the means of 
			pattern-counting statistics as such linear combinations. A
			procedure for calculating such expressions for statistics counting
			occurrences of classical and vincular patterns of length 3 is
			developed, and is then used to calculate all these expressions.
			The results can be used to 
			compute $\mathbb{E}_{\Gamma}(s,t)$ for all the above statistics, 
			and for all functions on $S_{n}$ that are linear combinations
			of them.
		
	\end{abstract}
	\maketitle

\section{Introduction}

Permutation patterns have been studied for decades, both in 
the field of computer science where Knuth \cite[Chapter 2.2.1]{knuth}
was one of the first to study avoidance of permutation patterns, and in
enumerative combinatorics (see e.g. \cite{simschm}). 
Over the years numerous articles have been written about permutation
patterns; mainly about permutations avoiding certain patterns.
An extensive treatment of patterns is given in \cite{kita}.
Random walks on Cayley graphs have been studied even
longer. They are used to model several different classes of problems.
One old problem is how many times a deck of cards must be 
shuffled with different methods until it is close to random. 
Another problem modelled by random walks is random number 
generation.
An introduction and overview of random walks on groups
is given by \cite{diac} and \cite{saloff}.

Here these subjects are combined by studying the values of 
different pattern-counting permutation statistics on the 
results of random 
walks on the Cayley graph of the symmetric group $\symg{n}$ 
induced by a subset $\Gamma$ of $\symg{n}$.
Such a walk with $t$ steps corresponds to a product of $t$ random 
elements of $\Gamma$. It is natural to be interested in the expected 
values of the permutation statistics on the product.

\begin{deff}
	Let $\Gamma$ be a subset of $\symg{n}$, $s$ a function from
	$\symg{n}$ to $\mathbb{R}$, and $t$ a non-negative integer. Then
	$\mathbb{E}_{\Gamma}(s,t)$ is the expected value of 
	$s(\gamma_{1} \cdots \gamma_{t})$, where every $\gamma_{i}$ is
	chosen independently from the uniform distribution on $\Gamma$.
\end{deff}

The behaviour of $\mathbb{E}_{\Gamma}(s,t)$ has been 
studied for permutation statistics like e.g.
the number of inversions, the absolute length function, and 
the fixed point counting function, usually for $\Gamma$ being some set
of transpositions. 
See for example \cite{bousquet}, \cite{eriksen}, \cite{erihu}, 
\cite{eriksson}, \cite{jonsson},
and \cite{sjostrand}. Expanding on \cite{erihu}, Hultman 
\cite{hult} found a method for computing $\mathbb{E}_{\Gamma}(s,t)$, 
given that $\Gamma$ is a union of conjugacy 
classes. He computed $\mathbb{E}_{\Gamma}(s,t)$ for some common permutation
statistics, such as $k$-cycle numbers, exceedances, inversion number,
descent number, and major index.

Even though both permutation patterns and random walks on
Cayley graphs have been studied from many different 
viewpoints, patterns have not yet been studied systematically 
in this context.
Therefore, permutation statistics counting occurrences of 
patterns are an interesting class to study.

This paper studies $\mathbb{E}_{\Gamma}(s,t)$ for 
permutation statistics that count the occurrences of 
different 3-patterns, using the method in \cite{hult}.
The 3-patterns investigated here are both the classical
3-patterns, and, more generally, the vincular 3-patterns. The 
vincular patterns can also be called generalised patterns, 
Babson-Steingr\'{i}msson patterns, or dashed patterns.
They were first described by Babson and Steingr\'{i}msson 
\cite{babste}. Steingr\'{i}msson \cite{steingr} has also 
written a survey of vincular patterns.

\begin{deff}
	An \emph{occurrence of a classical pattern} 
	$\phi = \phi_1 \phi_2 \ldots \phi_k \in \symg{k}$
	in a permutation $\pi = \pi_1 \pi_2 \ldots \pi_n \in \symg{n}$ 
	is a subsequence in $\pi$ of length $k$ whose letters are in the
	same relative order as those in $\phi$.
\end{deff}

For example, an occurrence of the classical pattern $123$ in 
$\pi \in \symg{n}$ is a subsequence $\pi_{i_{1}} \pi_{i_{2}} 
\pi_{i_{3}}$, where $i_{1}<i_{2}<i_{3}$, such that $\pi_{i_{1}}
< \pi_{i_{2}} < \pi_{i_{3}}$. 

\begin{deff}
	A \emph{vincular pattern} $\phi$ is written as a permutation 
	in $\symg{k}$ enclosed by brackets (square brackets or
	parentheses) which may have dashes 
	between adjacent letters. If adjacent letters are not
	separated by a dash, then the corresponding letters in an 
	occurrence of $\phi$ in $\pi \in \symg{n}$ must be adjacent.
	If $\phi$ begins with a square bracket then any occurrence
	of $\phi$ in $\pi$ must begin with $\pi_1$, and if $\phi$
	ends with a square bracket then any occurrence of $\phi$ in
	$\pi$ must end with $\pi_n$.
\end{deff}

For example, an occurrence of the vincular pattern $[1$-$23)$ 
in $\pi \in \symg{n}$ is a subsequence $\pi_{i_{1}} \pi_{i_{2}} 
\pi_{i_{3}}$, where $i_{1}<i_{2}<i_{3}$, $i_{1} = 1$, and
$i_{3} = i_{2}+1$, such that $\pi_{i_{1}}<\pi_{i_{2}}<\pi_{i_{3}}$. 

Note that the classical pattern $123$ can be written as the
vincular pattern \linebreak $(1$-$2$-$3)$. Hence the set of classical 
patterns is a subset of the set of vincular patterns.

The disposition of this article is as follows. Some theory
of the symmetric group and the method in \cite{hult} for computing 
$\mathbb{E}_{\Gamma}(s,t)$ are described in Section~\ref{theory}.
The only prerequisite for being able to use that method is to express
the mean statistic corresponding to $s$ as a linear combination 
of irreducible characters of $\symg{n}$.
The main focus of this article is to calculate such expressions.
A procedure for this is developed in 
Sections~\ref{chicoeffs} and \ref{proced}. In 
Section~\ref{chicoeffs} some notation is introduced together
with a couple of useful lemmas. Then in Section~\ref{proced}
the actual procedure is described and proved. It is used to 
calculate the coefficients for all vincular 3-patterns. 
It turns out that only five irreducible characters are needed.

The procedure should be fairly easy to generalise to longer
patterns, even though that would result in more computational 
work. An interesting question for future studies is how many 
irreducible characters that are needed in the expressions 
for means of statistics counting vincular $k$-patterns. 
It seems that the only irreducible characters needed are 
those that are indexed by partitions 
$\lambda=(\lambda_1, \lambda_2,\ldots,\lambda_{\ell})$, where
$\lambda_1 \geq n-k$ and $\lambda_2 < k$.

In Section~\ref{expvals}, the expected values are computed
explicitly for the statistics counting classical 3-patterns 
when $\Gamma$ is the set of transpositions.

\section{Preliminaries} \label{theory}

First some theory of the symmetric group is described.
For a more thorough treatment, see e.g. \cite{sagan}.

Let $\symg{n}$ denote the symmetric group of permutations
of $[n]=\{1, \ldots,n\}$.
The sequence $\lambda = (\lambda_1,\lambda_2,\ldots,\lambda_{\ell})$
of the lengths of all cycles in $\pi$
in weakly decreasing order is called the cycle type
of $\pi$. The sequence $\lambda$ is an integer partition of $n$, 
denoted $\lambda \vdash n$.
The set of integer partitions of $n$ will be called $P_{n}$.

Each conjugacy class in $\symg{n}$ consists of all permutations
of a given cycle type. The conjugacy class with cycle type
$\lambda \vdash n$ is denoted
$C_{\lambda}$.
Functions from $P_{n}$ to $\mathbb{C}$ are called class
functions. Class functions can also be seen as functions 
from $\symg{n}$ to $\mathbb{C}$ which are constant on 
conjugacy classes. 

There is a standard bijection between $P_{n}$ and the irreducible 
representations of $\symg{n}$. The character corresponding 
to the representation indexed by $\lambda \vdash n$ will be
denoted $\chi^{\lambda}$. These irreducible characters form an
orthonormal basis for the space of class functions with respect 
to the standard Hermitian inner product
\begin{displaymath}
	\langle f,g \rangle = \frac{1}{n!} \sum_{\lambda \vdash n}
	|C_{\lambda}| f(\lambda) g(\lambda)^{*}, 
\end{displaymath}
where $g(\lambda)^{*}$ is the complex conjugate of $g(\lambda)$.

The procedure described in \cite{hult} for computing 
$\mathbb{E}_{\Gamma}(s,t)$ will be summarised here.

\begin{deff}
	Take a function $s: \symg{n} \to \mathbb{R}$.
	The \emph{mean statistic} $\overline{s}$ is the class
	function computing the mean of $s$ over conjugacy classes, i.e.
	\begin{displaymath}
		\overline{s}(\lambda) = 
		\frac{1}{|C_{\lambda}|} \sum_{\pi \in C_{\lambda}} s(\pi).
	\end{displaymath}
	Note that if $s$ is a class function, then $s = \overline{s}$.
\end{deff}

Let $n_{t}$ be the permutation statistic defined by
\begin{displaymath}
	n_{t}(\pi) = \bigl\lvert \{(\gamma_{1}, \ldots,\gamma_{t}) \in \Gamma^{t}|
	\gamma_{1} \cdots \gamma_{t} = \pi \} \bigr\rvert .
\end{displaymath}
This is a class function for all $t$ if and only if
$\Gamma$ is a union of conjugacy classes.

The following theorem is proved in \cite{hult}.
\begin{thm} \label{hultone}
	If at least one of the statistics $n_{t}$ and $s$ is a class 
	function, then
	\begin{displaymath}
		\mathbb{E}_{\Gamma}(s,t) = \frac{n!}{|\Gamma|^{t}} \langle
		\overline{s}, \overline{n_{t}} \rangle.
	\end{displaymath}
\end{thm}

The inner product $\langle \overline{s}, \overline{n_{t}} 
\rangle$ is easy to compute if $\overline{s}$ and 
$\overline{n_{t}}$ are expressed in the basis of irreducible
$\symg{n}$-characters. Let the coefficients $a_{\lambda}$ and
$b_{\lambda}^{(t)}$ be defined by
\begin{align*}
	\overline{s} & = \sum_{\lambda \vdash n} a_{\lambda} 
	\chi^{\lambda}, &
	\overline{n_{t}} & = \sum_{\lambda \vdash n} b_{\lambda}^{(t)}
	\chi^{\lambda}.
\end{align*}
Then
\begin{displaymath}
	\mathbb{E}_{\Gamma}(s,t) = \frac{n!}{|\Gamma|^{t}}
	\sum_{\lambda \vdash n} a_{\lambda} b_{\lambda}^{(t)}.
\end{displaymath}

The following theorem can be used to compute the coefficients
$b_{\lambda}^{(t)}$ when $\Gamma$ is a conjugacy class. It follows
easily from \cite[Theorem A.1.9]{lanzvo}. There is a similar formula
for $b_{\lambda}^{(t)}$ when $\Gamma$ is a union of conjugacy 
classes which can be found e.g. in \cite{hult}.

\begin{thm} \label{hulttwo}
	Suppose that $\Gamma$ is a conjugacy class. Let
	$\mu$ denote the cycle type of the permutations in $\Gamma$.
	Then
	\begin{displaymath}
		b_{\lambda}^{(t)} = \frac{1}{n!\Bigl(\chi^{\lambda}\bigl((1,\ldots,1)
			\bigr)\Bigr)^{t-1}} \Bigl(|\Gamma|
		\chi^{\lambda}(\mu) \Bigr)^{t}.
	\end{displaymath}
\end{thm}

This means that if the coefficients $a_{\lambda}$ in the 
expression $\overline{s}  = \sum_{\lambda \vdash n} a_{\lambda}
\chi^{\lambda}$ are found, then $\mathbb{E}_{\Gamma}(s,t)$ can be 
computed for every $\Gamma$ being a union of conjugacy 
classes.
Therefore, the main focus of this article is a systematic
method for calculating such expressions.
The next section shows how to do that for a great number of pattern
statistics.

\section{Mean statistics of 3-patterns} \label{chicoeffs}

In this section and the next the means of permutation 
statistics counting the number of vincular 3-patterns 
(including classical 3-patterns) are investigated. A
procedure for finding expressions for the mean statistics
as linear combinations of the irreducible characters of
$\symg{n}$ is described. The procedure is used to 
find such expressions for the mean statistics of all 
vincular 3-patterns. 

Let $\patt{\phi}(\pi)$ denote the number of occurrences of
the pattern $\phi$ in $\pi$. 

If $s_{1}$ and $s_{2}$ are permutation statistics, then it
follows from the definition that $\overline{\alpha s_{1}+
	\beta s_{2}} = \alpha \overline{s_{1}}+ \beta \overline{s_{2}}$,
where $\alpha$ and $\beta$ are constants. 
Hence the results of this section and the next can be 
applied in order to
get the expressions for the means of all
functions on $\symg{n}$ that are linear combinations of
the permutation statistics counting vincular
3-patterns. 
As an illustration, consider the following statistic:
A peak of $\pi \in \symg{n}$ 
is an index $i$ such that $\pi(i-1) < \pi(i) > \pi(i+1)$.
Let $\peak(\pi)$ denote 
the number of peaks of $\pi$. 
Then $\peak(\pi) = \patt{(132)}(\pi)+
\patt{(231)}(\pi)$, and it follows that $\overline{\peak} = 
\overline{\patt{(132)}} + \overline{\patt{(231)}}$.

Babson and Steingr\'{i}msson \cite{babste} showed that
many of the known Mahonian statistics can be written as sums of
three permutation statistics counting vincular 3-patterns and
a permutation statistic counting descents ($\des$) or
ascents ($\asc$). 
It is known that $\overline{\des} = 
\frac{n-1}{2} \chi^{(n)} - \frac{1}{n} \chi^{(n-1,1)} - 
\frac{1}{n} \chi^{(n-2,1,1)}$; see \cite[Theorem 2.1]{fulman} 
or \cite{hult}. It easily follows that $\overline{\asc} = 
\frac{n-1}{2} \chi^{(n)} + \frac{1}{n} \chi^{(n-1,1)} + 
\frac{1}{n} \chi^{(n-2,1,1)}$.
Hence, the means of those Mahonian statistics 
can be analysed using the results of the present paper.

This section contains some definitions and lemmas, and in the
next section the procedure mentioned above is described. From now
on only permutations of $[n]$, where $n \geq 3$, will be
considered, since shorter permutations cannot contain
3-patterns.

Let $\phi$ be a vincular 3-pattern.
As will be seen in Section~\ref{proced}, for all permutation 
statistics $\patt{\phi}$, the means
$\overline{\patt{\phi}}$ are linear combinations 
of five of the seven irreducible characters appearing in the following 
lemma. The remaining two characters are necessary for the proofs, 
but disappear in the results.

\begin{lem} \label{chin}
	Let $\lambda \vdash n$. If $\lambda$ has $p$ parts of size 1, 
	$q$ parts of size 2, and $r$ parts of size 3, then 
	\begin{align*}
		\chi^{(n)}(\lambda) & = 1,\\
		\chi^{(n-1,1)}(\lambda) & = p-1,\\
		\chi^{(n-2,2)}(\lambda) & = \binom{p-1}{2}+q-1,\\
		\chi^{(n-2,1,1)}(\lambda) & = \binom{p-1}{2}-q,\\
		\chi^{(n-3,3)}(\lambda) & = \binom{p-1}{3}+(p-1)(q-1)+r,\\
		\chi^{(n-3,2,1)}(\lambda) & = 2\binom{p-1}{3}-p-r+2,\\
		\chi^{(n-3,1,1,1)}(\lambda) & = \binom{p-1}{3}-q(p-1)+r.
	\end{align*}
\end{lem}

\begin{bevis}
	This easily follows from the Murnaghan-Nakayama rule (see e.g.
	\cite[Chapter 7.17]{stanley2}).
\end{bevis}

Note that if $n\leq 5$, then some of the functions above
are not irreducible characters of $\symg{n}$. In that case
they can be expressed as linear combinations of the other
functions above as follows: 
\begin{itemize}
	\item If $n=5$ then 
	$\chi^{(n-3,3)}(\lambda) = 0$.
	\item If $n=4$ then 
	$\chi^{(n-3,3)}(\lambda) = - \chi^{(n-2,2)}(\lambda)$ and
	$\chi^{(n-3,2,1)}(\lambda)=0$.
	\item If $n=3$ then 
	$\chi^{(n-2,2)}(\lambda)=0$, \ $\chi^{(n-3,3)}(\lambda) = 
	-\chi^{(n-1,1)}(\lambda)$, \newline $\chi^{(n-3,2,1)}(\lambda) =
	-\chi^{(n-2,1,1)}(\lambda)$, and $\chi^{(n-3,1,1,1)}(\lambda)=0$.
\end{itemize}

It is convenient to introduce the column vector
\begin{displaymath}
	\chivek(\lambda) = \begin{pmatrix} 
		\chi^{(n)}(\lambda) \\ \chi^{(n-1,1)}(\lambda) \\ \chi^{(n-2,2)}(\lambda)
		\\ \chi^{(n-2,1,1)}(\lambda) \\ \chi^{(n-3,3)}(\lambda) \\ 
		\chi^{(n-3,2,1)}(\lambda) \\ \chi^{(n-3,1,1,1)}(\lambda)
	\end{pmatrix}.
\end{displaymath}

To find the mean $\overline{\patt{\phi}}(\lambda)$, 
the total number of occurrences of $\phi$ in all permutations
$\pi \in C_{\lambda}$ will be calculated and the result will be
divided by $|C_{\lambda}|$.
For $1\leq i_{1} < i_{2} < i_{3}
\leq n$, define the indicator variable
\begin{displaymath}
	\patt{\phi}^{i_{1}, i_{2}, i_{3}}(\pi) = \begin{cases}
		1 & \text{if } \pi_{i_{1}} \pi_{i_{2}} \pi_{i_{3}} \text{ is an
			occurrence of } \phi, \\
		0 & \text{otherwise.}
	\end{cases}
\end{displaymath}
Thus $\patt{\phi}(\pi) = \sum_{i_{1}, i_{2}, i_{3}} 
\patt{\phi}^{i_{1}, i_{2}, i_{3}}(\pi)$.

Let $\fourset{j}(i_{1}, i_{2}, i_{3}) \subseteq C_{\lambda}$
be the set of all permutations $\pi \in C_{\lambda}$
such that $\{\pi_{i_{1}},\pi_{i_{2}},\pi_{i_{3}}\} \cup 
\{i_{1},i_{2},i_{3}\}$ has cardinality $j$. Then, 
for each choice of $i_{1}$, $i_{2}$, and $i_{3}$, 
the set $C_{\lambda}$ is the disjoint union of
$\fourset{j}(i_{1}, i_{2}, i_{3})$, $j=3,4,5,6$.

As will be seen, it is convenient to calculate the total
number of occurrences of $\phi$ in $C_{\lambda}$ as
\begin{displaymath}
	\sum_{\pi\in C_{\lambda}} \patt{\phi}(\pi) = 
	\sum_{j=3}^{6}  \Biggl(\sum_{i_{1}, i_{2}, i_{3}} \biggl(\sum_{\pi \in 
		\fourset{j}(i_{1}, i_{2}, i_{3})} 
	\patt{\phi}^{i_{1}, i_{2}, i_{3}}(\pi)\biggr)\Biggr).
\end{displaymath}

\begin{deff}
	For fixed indices $i_{1}<i_{2}<i_{3}$, and for a fixed 
	$\lambda \vdash n$, let an \emph{$M$-set} $M(x,y,z)$ be 
	the set of permutations 
	$\{\pi \in C_{\lambda} \ | \ \pi_{i_{1}}=x, \ \pi_{i_{2}}=y,
	\ \pi_{i_{3}}=z\}$.
\end{deff}

It is easy to see that each set 
$\fourset{j}(i_{1}, i_{2}, i_{3})$ is a (disjoint) union of 
$M$-sets. 

Note that for all permutations $\pi$ in an $M$-set, 
the subsequences $\pi_{i_{1}} \pi_{i_{2}} \pi_{i_{3}}$ are identical
and are occurrences of one of the six classical 3-patterns. 
Therefore the \emph{pattern of the $M$-set} will be defined as
the (classical) pattern given by $\pi_{i_{1}} \pi_{i_{2}} \pi_{i_{3}}$. 

The sum of the cardinalities of all $M$-sets in 
$\fourset{j}(i_{1}, i_{2}, i_{3})$ with the 
pattern given by the letters in $\phi$  is equal to
$\sum_{\pi \in 
	\fourset{j}(i_{1}, i_{2}, i_{3})} 
\patt{\phi}^{i_{1}, i_{2}, i_{3}}(\pi)$ if $i_{1}$, $i_{2}$, and 
$i_{3}$ satisfy the conditions given by the dashes and square brackets 
in $\phi$. The following lemma gives the cardinalities of all 
$M$-sets.

Let $\fall{n}{\ell}$ denote the falling factorial 
$n(n-1) \cdots (n-\ell+1)$.
\pagebreak

\begin{lem} \label{mcard}
	Suppose that indices $i_{1}<i_{2}<i_{3}$ and $\lambda \vdash n$ are given, 
	where $\lambda$ has $p$ parts of size 1, 
	$q$ parts of size 2, and $r$ parts of size 3. 
	Then the various $M$-sets fall into ten equicardinal symmetry 
	classes. The classes and their cardinalities are as follows,
	where $a_{4},a_{5},a_{6} \notin \{i_{1},i_{2},i_{3}\}$:
	\begin{itemize}
		\item[1)] $|M(i_{1},i_{2},i_{3})| = H_{31}(\lambda) 
		|C_{\lambda}|/\fall{n}{3}$,
		\item[2)] $|M(i_{1},i_{3},i_{2})| = |M(i_{3},i_{2},i_{1})| = 
		|M(i_{2},i_{1},i_{3})| = H_{32}(\lambda)|C_{\lambda}|/
		\fall{n}{3}$,
		\item[3)] $|M(i_{2},i_{3},i_{1})| = |M(i_{3},i_{1},i_{2})| =  
		H_{33}(\lambda) |C_{\lambda}|/\fall{n}{3}$,
		\item[4)] $|M(i_{1},i_{2},a_{4})| = |M(i_{1},a_{4},i_{3})| = |M(a_{4},i_{2},i_{3})| = 
		H_{41}(\lambda) |C_{\lambda}|/\fall{n}{4}$,
		\item[5)] $|M(i_{2},i_{1},a_{4})| = |M(i_{3},a_{4},i_{1})| = |M(a_{4},i_{3},i_{2})| = 
		H_{42}(\lambda) |C_{\lambda}|/\fall{n}{4}$,
		\item[6)] $|M(i_{1},i_{3},a_{4})| = |M(i_{1},a_{4},i_{2})| = |M(i_{3},i_{2},a_{4})| = 
		|M(a_{4},i_{2},i_{1})| = \newline |M(i_{2},a_{4},i_{3})| = |M(a_{4},i_{1},i_{3})| = 
		H_{43}(\lambda) |C_{\lambda}|/\fall{n}{4}$,
		\item[7)] $|M(i_{2},i_{3},a_{4})| = |M(i_{3},a_{4},i_{2})| = |M(i_{3},i_{1},a_{4})| = 
		|M(a_{4},i_{3},i_{1})| = \newline |M(i_{2},a_{4},i_{1})| = |M(a_{4},i_{1},i_{2})| = 
		H_{44}(\lambda) |C_{\lambda}|/\fall{n}{4}$,
		\item[8)] $|M(i_{1},a_{4},a_{5})| = |M(a_{4},i_{2},a_{5})| = |M(a_{4},a_{5},i_{3})|  
		= H_{51}(\lambda) |C_{\lambda}|/\fall{n}{5}$,
		\item[9)] $|M(i_{2},a_{4},a_{5})| = |M(i_{3},a_{4},a_{5})| = |M(a_{4},i_{1},a_{5})| = 
		|M(a_{4},i_{3},a_{5})| = \newline |M(a_{4},a_{5},i_{1})| = |M(a_{4},a_{5},i_{2})| = 
		H_{52}(\lambda) |C_{\lambda}|/\fall{n}{5}$,
		\item[10)] $|M(a_{4},a_{5},a_{6})| = H_{61}(\lambda) 
		|C_{\lambda}|/\fall{n}{6}$,
	\end{itemize}
	where each $H_{j\ell}(\lambda) |C_{\lambda}|/\fall{n}{j}$ with $j>n$ is
	to be replaced by $0$, and $H_{j\ell}(\lambda)$ are given by
	\begin{itemize}
		\item[1)] $H_{31}(\lambda) = p(p-1)(p-2)$,
		\item[2)] $H_{32}(\lambda) = 2pq$,
		\item[3)] $H_{33}(\lambda) = 3r$,
		\item[4)] $H_{41}(\lambda) = p(p-1)(n-p)$,
		\item[5)] $H_{42}(\lambda) = 2q(n-p-2)$,
		\item[6)] $H_{43}(\lambda) = p(n-p-2q)$,
		\item[7)] $H_{44}(\lambda) = n-p-2q-3r$,
		\item[8)] $H_{51}(\lambda) = p^3+3p^2+2pq-n(2p^2+3p)+n^2p$,
		\item[9)] $H_{52}(\lambda) = p^2+4p+2pq+8q+3r-n(2p+2q+4)+n^2$,
		\item[10)] $H_{61}(\lambda) = -p^3-9p^2-20p-6pq-24q-6r
		+n(3p^2+18p+6q+20)$ \newline $\mbox{ }\qquad\qquad -n^2 (3p+9)+n^3$.
	\end{itemize}
	
	The expressions defining $H_{j\ell}(\lambda)$ are $0$ when $j>n$.
\end{lem}

\begin{bevis}
	
	Writing the permutations of the $M$-sets in cycle notation 
	easily shows that they have the same 
	cardinality. Hence it is enough to show that one of the 
	$M$-sets in each class has the stated cardinality.
	
	Cases (2), (4), (8), and (10) will be shown. The rest have similar 
	proofs which will be omitted here.
	
	Let $\tau_{\ell}$ be the number of $\ell$-parts in $\lambda$ for
	$4 \leq \ell \leq n$.
	Recall that 
	\begin{displaymath}
		|C_{\lambda}| = \frac{n!}{1^p p! 2^q q! 3^r r! 
			4^{\tau_{4}} \tau_{4}! \cdots n^{\tau_{n}} \tau_{n}!}.
	\end{displaymath}
	
	\noindent {\bf Case (2)}: The permutations
	in $M(i_{1},i_{3},i_{2})$ written in cycle notation look like
	$(i_{1})(i_{2} \; i_{3}) \cdots$. There are as many ways to 
	complete the permutations as the number of permutations in the
	conjugacy class $C_{\mu}$ of $\symg{n-3}$, where $\mu$ has
	$(p-1)$ $1$-parts, $(q-1)$ $2$-parts, $r$ $3$-parts, and $\tau_{\ell}$
	$\ell$-parts for all $\ell \geq 4$.
	Hence 
	\begin{displaymath}
		|M(i_{1},i_{3},i_{2})| = |C_{\mu}| = \frac{(n-3)!}{1^{p-1} (p-1)! 2^{q-1} 
			(q-1)! 3^r r! 4^{\tau_{4}} \tau_{4}! \cdots n^{\tau_{n}} \tau_{n}!} = 
		\frac{2pq}{\fall{n}{3}} |C_{\lambda}|.
	\end{displaymath}
	
	If $p=0$ or $q=0$, then $C_{\mu}$ does not exist and 
	$|M(i_{1},i_{3},i_{2})|=0 = \frac{2pq}{\fall{n}{3}} |C_{\lambda}|$.
	
	\vspace{2mm}
	\noindent Cases (1) and (3) can be proved in the same way as Case (2).
	
	\vspace{2mm}
	\noindent {\bf Case (4)}: The permutations in $M(i_{1},i_{2},a_{4})$ are exactly 
	all permutations in $C_{\lambda}$ that can be written as 
	$(i_{1})(i_{2})(i_{3} \; a_{4} \;\ldots) \cdots$ in cycle notation.
	
	Let $\|(i_{1}) (i_{2})(i_{3} \; a_{4} \; \ldots)\cdots\|_{\lambda}$ 
	denote the number of permutations in $C_{\lambda}$ that can be 
	written as $(i_{1}) (i_{2})(i_{3} \; a_{4} \; \ldots)\cdots$,
	etc. 
	
	It is easy to see that the number of permutations $\pi$ in $C_{\lambda}$
	where $\pi(i_{1}) = i_{1}$, $\pi(i_{2}) = i_{2}$, and 
	$\pi(i_{3}) \neq i_{3}$ is equal to 
	$ \|(i_{1}) (i_{2})\cdots\|_{\lambda} - 
	\|(i_{1}) (i_{2})(i_{3})\cdots\|_{\lambda}.$
	This is of course $0$ if $n=3$.
	There are $(n-3)$ possible values of $\pi(i_{3})$ if $n\geq 4$, 
	and all of them are equally likely. Hence 
	\begin{displaymath}
		|M(i_{1},i_{2},a_{4})| = \frac{1}{n-3}
		\bigl(\|(i_{1}) (i_{2})\cdots\|_{\lambda} - 
		\|(i_{1}) (i_{2})(i_{3})\cdots\|_{\lambda}\bigr).
	\end{displaymath}
	
	Similar calculations as in Case (2) give that
	\begin{align*}
		\|(i_{1}) (i_{2})\cdots\|_{\lambda} - 
		\|(i_{1}) (i_{2})(i_{3})\cdots\|_{\lambda} & = 
		\frac{p(p-1)}{\fall{n}{2}}|C_{\lambda}| - 
		\frac{p(p-1)(p-2)}{\fall{n}{3}}|C_{\lambda}| \\
		& =
		\frac{p(p-1)(n-p)}{\fall{n}{3}}|C_{\lambda}|.
	\end{align*}

	Thus, if $n\geq 4$ then $|M(i_{1},i_{2},a_{4})| =p(p-1)(n-p)
	|C_{\lambda}|/\fall{n}{4}$, and if $n=3$ then $|M(i_{1},i_{2},a_{4})|=0$
	and $p(p-1)(n-p)=0$.
	
	\vspace{2mm}
	\noindent Cases (5) -- (7) can be shown in the same way as Case (4).

	\vspace{2mm}
	\noindent {\bf Case (8)}: The permutations in $M(i_{1},a_{4},a_{5})$ are exactly
	all permutations in $C_{\lambda}$ that can be written as 
	$(i_{1})(i_{2} \; a_{4} \ldots i_{3} \;a_{5} \ldots) \cdots$
	or $(i_{1})(i_{2} \;a_{4} \ldots)(i_{3} \; a_{5} \ldots) \cdots$.
	Of course there are no such permutations if $n<5$.
	
	Reasoning as in Case (4) when $n \geq 5$, it is straightforward to get
	\begin{align*}
		|M(i_{1},a_{4},a_{5})| & = \frac{1}{\fall{(n-3)}{2}} 
		\biggl(\|(i_{1})\cdots\|_{\lambda} -
		\|(i_{1})(i_{2})(i_{3})\cdots\|_{\lambda} -
		\|(i_{1})(i_{2} \; i_{3})\cdots\|_{\lambda} \\
		& \quad -
		(n-3)\Bigl(\|(i_{1})(i_{2})(i_{3} \; \alpha \; \ldots)\cdots\|_{\lambda}+
		\|(i_{1})(i_{3})(i_{2} \; \alpha \; \ldots)\cdots\|_{\lambda} \\
		& \quad +
		\|(i_{1})(i_{2} \; i_{3} \; \alpha \; \ldots)\cdots\|_{\lambda}+
		\|(i_{1})(i_{3} \; i_{2} \; \alpha \; \ldots)\cdots\|_{\lambda}\Bigr)\biggr),
	\end{align*}
	where the quantity inside the parenthesis is equal to $0$ if
	$n=3$ or $n=4$. 
	
	It easily follows that if $n\geq 5$ then
	
	\begin{align*}
		|M(i_{1},a_{4},a_{5})| & = \frac{1}{\fall{(n-3)}{2}} \Bigl(
		\frac{p}{n}|C_{\lambda}|-|M(i_{1},i_{2},i_{3})| - |M(i_{1},i_{3},i_{2})|\\
		& \quad - 
		(n-3)\bigl( 2 |M(i_{1},i_{2},\alpha)|+2|M(i_{1},i_{3},\alpha)| \bigr)
		\Bigr)\\
		& = (p^3+3p^2+2pq-n(2p^2+3p)+n^2p) \frac{|C_{\lambda}|}{\fall{n}{5}},\;
	\end{align*}
	and if $n=3$ or $n=4$ then $(p^3+3p^2+2pq-n(2p^2+3p)+n^2p)=0$ and 
	$|M(i_{1},a_{4},a_{5})|=0$.
	
	\vspace{2mm}
	\noindent Case (9) can be shown in the same way as Case (8).
	
	\vspace{2mm}
	\noindent {\bf Case (10)}:
	The conjugacy class $C_{\lambda}$ is the union of all 
	$M$-sets. There are $\fall{(n-3)}{3}$ $M$-sets in the
	last class, and it is obvious that they have the same 
	cardinality.
	Hence, if $n\geq 6$,
	\begin{align*}
		|M(a_{4},a_{5},a_{6})| & = \frac{1}{\fall{(n-3)}{3}}|C_{\lambda}|
		\biggl(1 - \frac{H_{31}(\lambda)}{\fall{n}{3}}- 
		3 \frac{H_{32}(\lambda)}{\fall{n}{3}}- 
		2 \frac{H_{33}(\lambda)}{\fall{n}{3}}\\
		& \quad - (n-3) \Bigl( 
		3\frac{H_{41}(\lambda)}{\fall{n}{4}} 
		+3\frac{H_{42}(\lambda)}{\fall{n}{4}}
		+6\frac{H_{43}(\lambda)}{\fall{n}{4}} 
		+6\frac{H_{44}(\lambda)}{\fall{n}{4}}\Bigr) \\
		& \quad - (n-3)(n-4) \Bigl(3\frac{H_{51}(\lambda)}{\fall{n}{5}}
		+ 6\frac{H_{52}(\lambda)}{\fall{n}{5}}\Bigr)\biggr)\\
		& = \bigl( -p^3-9p^2-20p-6pq-24q-6r
		+n(3p^2+18p+6q+20)\\
		& \quad -n^2 (3p+9)+n^3 \bigr) \frac{|C_{\lambda}|}{\fall{n}{6}},
	\end{align*}
	and if $n \leq 5$ then $|M(a_{4},a_{5},a_{6})|=0$.
	
	It is easy to see that the expression inside the first (and
	therefore the last)
	parenthesis is $0$ if $n=5$.
	Since $H_{j\ell}(\lambda) = 0$ if $5\geq j > n$ it also follows
	that the expression in the last parenthesis is $0$ if
	$n=3$ or $n=4$. Hence $H_{61}(\lambda) = 0$ for all $n \leq 5$.
	
	Thus all $M$-sets have the cardinalities 
	stated, and $H_{j\ell} = 0$ if $j>n$.
\end{bevis}

\begin{lem} \label{hlincomb}
	The functions $H_{j\ell}(\lambda)$ defined in Lemma~\ref{mcard} 
	are linear combinations of the seven irreducible characters in 
	Lemma~\ref{chin}, where the coefficients are as follows:
	\begin{align*} 
		H_{31}(\lambda) &  = (
		1 \vecsp 3 \vecsp 3 \vecsp 3 \vecsp 1 \vecsp 2 \vecsp 1
		) \, \chivek(\lambda),\\ 
		H_{32}(\lambda) & = (
		1 \vecsp 1 \vecsp 1 \vecsp -1 \vecsp 1 \vecsp 0 \vecsp -1
		) \, \chivek(\lambda), \\ 
		H_{33}(\lambda) &  = (
		1 \vecsp 0 \vecsp 0 \vecsp 0 \vecsp 1 \vecsp -1 \vecsp 1
		) \, \chivek(\lambda),\\
		H_{41}(\lambda) &= (
		n-3 \vecsp 2n-7 \vecsp n-5 \vecsp n-5 \vecsp -1 \vecsp -2 \vecsp -1
		) \, \chivek(\lambda),\\
		H_{42}(\lambda) &= (
		n-3 \vecsp -1 \vecsp n-3 \vecsp -n+3 \vecsp -1 \vecsp 0 \vecsp 1
		) \, \chivek(\lambda),\\
		H_{43}(\lambda) &= (
		n-3 \vecsp n-4 \vecsp -2 \vecsp 0 \vecsp -1 \vecsp 0 \vecsp 1
		) \, \chivek(\lambda), \\ 
		H_{44}(\lambda) &= (
		n-3 \vecsp -1 \vecsp -1 \vecsp 1 \vecsp -1 \vecsp 1 \vecsp -1
		) \, \chivek(\lambda),\\
		H_{51}(\lambda) &= (
		n^2-7n+12 \vecsp n^2-9n+20 \vecsp -2n+10 \vecsp -2n+8 \vecsp 2 
		\vecsp 2 \vecsp 0) \, \chivek(\lambda), \\ 
		H_{52}(\lambda) &=(
		n^2-7n+12 \vecsp -2n+8 \vecsp -n+6 \vecsp n-4 \vecsp 2 \vecsp -1 \vecsp 0
		) \, \chivek(\lambda),\\
		H_{61}(\lambda) & =(
		n^3-12 n^2+47n-60 \vecsp -3 n^2+27n-60 \vecsp 6n-30 \vecsp 0 \vecsp -6 \vecsp 0 \vecsp 0
		) \, \chivek(\lambda).
	\end{align*}
	
	If $n \leq 5$, these expressions can be transformed to 
	linear combinations of the irreducible characters of 
	$\symg{n}$, as described after Lemma~\ref{chin}.
\end{lem}

\begin{bevis}
	Follows directly from Lemmas~\ref{chin} and~\ref{mcard}.
\end{bevis}

Since the permutations in each $M$-set have the same
letters at positions $i_{1}$, $i_{2}$, and $i_{3}$, it
follows from Lemma~\ref{mcard} that
the total number of occurrences of $\phi$ in $C_{\lambda}$ 
is a sum where all terms are of the form 
$H_{j \ell}(\lambda) |C_{\lambda}|/ \fall{n}{j}$.
Dividing by $|C_{\lambda}|$ gives the mean
$\overline{\patt{\phi}}$ as a linear combination of the 
functions $H_{j \ell}(\lambda)$. Using Lemma~\ref{hlincomb},
$\overline{\patt{\phi}}$ can then be written as a
linear combination of irreducible characters.

\section{A procedure for calculating the means of vincular 3-patterns}
\label{proced}

This section describes a procedure to calculate the
means $\overline{\patt{\phi}}$ as linear combinations
of the irreducible characters in Lemma~\ref{chin}. 

The procedure will then be used to calculate the
coefficients for all (classical and) vincular 3-patterns.

Recall that 
\begin{displaymath}
	\overline{\patt{\phi}}(\pi) = \frac{1}{|C_{\lambda}|}
	\sum_{\pi \in C_{\lambda}} \patt{\phi}(\pi).
\end{displaymath}
This sum will be computed as
\begin{displaymath}
	\frac{1}{|C_{\lambda}|} \sum_{j=3}^{6}  \Biggl(\sum_{i_{1}, i_{2}, i_{3}} 
	\biggl(\sum_{\pi \in \fourset{j}(i_{1}, i_{2}, i_{3})} 
	\patt{\phi}^{i_{1}, i_{2}, i_{3}}(\pi)\biggr) \Biggr).
\end{displaymath}

First some notation will be introduced. 
\begin{itemize}
	\item Let $R_{\lambda,j}(\phi)$ denote the sum 
	$\sum_{i_{1}, i_{2}, i_{3}} \sum_{\pi \in 
		\fourset{j}(i_{1}, i_{2}, i_{3})} 
	\patt{\phi}^{i_{1}, i_{2}, i_{3}}(\pi)$.
	
	\item The classical pattern obtained from the
	vincular pattern $\phi$ by removing all brackets and dashes will 
	be denoted $\phi_{c}$. 
	
	\item The conditions on $i_{1}$, $i_{2}$, 
	and $i_{3}$ given by the square brackets and absence of dashes in 
	$\phi$ together with the condition
	$i_{1}<i_{2}<i_{3}$ will be called $CI_{\phi}$. 
	
	\item Let $k_{\phi}$ denote the number of conditions
	given by the square brackets and absence of dashes in $\phi$. 
	It will be assumed that $k_{\phi}\leq 3$. This will not
	result in any loss of information, since patterns of the
	form $[abc]$ only occur when $n=3$, and in that case 
	$[abc]$ is identical to $(abc)$.
	
\end{itemize}

As an example, if $\phi$ is the vincular pattern $[\mbox{$21$-$3$})$, 
then $\phi_{c}$ is the classical pattern $213$ and
$CI_{\phi}$ is comprised of the conditions 
$i_{1}<i_{2}<i_{3}$, $i_{1} = 1$, and 
$i_{2} = i_{1}+1$.

Remember that $\sum_{\pi \in \fourset{j}(i_{1}, i_{2}, i_{3})} 
\patt{\phi}^{i_{1}, i_{2}, i_{3}}(\pi)$ is equal to the sum of
the cardinalities of all $M$-sets in 
$\fourset{j}(i_{1}, i_{2}, i_{3})$ with pattern $\phi_{c}$.
Summing over all $i_{1}$, $i_{2}$, and $i_{3}$
satisfying $CI_{\phi}$ gives $R_{\lambda,j}(\phi)$.

In order to facilitate the calculation of 
$R_{\lambda,j}(\phi)$, the following function is defined.

\begin{deff}
	Let $\indfn_{\ell}(\phi)$ be the number of index suites
	$(i_{1}, i_{2}, i_{3}, y_{4}, \ldots, y_{\ell})$ satisfying
	$CI_{\phi}$ and a condition stating the relative order of the 
	indices (an order not violating the conditions $CI_{\phi}$).
\end{deff}

If for example, the pattern $\phi = [\mbox{$21$-$3$})$ and
$\ell = 5$, then $\indfn_{5}(\phi)$ is the number of 
index suites $(i_{1}, i_{2}, i_{3}, y_{4}, y_{5})$ satisfying
$i_{1}<i_{2}<i_{3}$, $i_{1} = 1$, $i_{2} = i_{1}+1$, and (e.g.)
$i_{1}<i_{2}<y_{5}<i_{3}<y_{4}$. An example of a relative 
order of the indices violating $CI_{\phi}$ is 
$i_{1}<y_{4}<i_{2}<y_{5}<i_{3}$, since the condition 
$i_{2} = i_{1}+1$ is included in $CI_{\phi}$.

\begin{lem} \label{indf}
	The function $\indfn_{\ell}(\phi)$ is well defined, and
	is given by
	\begin{displaymath}
		\indfn_{\ell}(\phi) = \binom{n-k_{\phi}}{\ell-k_{\phi}}.
	\end{displaymath}
\end{lem}

\begin{bevis}
	To choose $\ell$ indices with a given relative order 
	in $[n]$ is equivalent to choose $\ell+1$
	spacings between the indices (including the spaces before
	the first index and after the last index) whose sum
	is $n-\ell$.
	Every condition given by square brackets or absence of dashes in 
	$\phi$ decreases the number of spacings to choose by 1, but
	their sum should still be $n-\ell$. 
	Hence $\indfn_{\ell}(\phi)$ is the number of ways
	to put $n-\ell$ identical objects in $\ell+1-k_{\phi}$
	different boxes (spacings), that is 
	$\ttbinom{n-k_{\phi}}{\ell-k_{\phi}}$.
\end{bevis}

Now it will be shown how to calculate $R_{\lambda,j}(\phi)$
for each $j$. Note that the mean $\overline{\patt{\phi}} =  
\bigl(R_{\lambda,3}(\phi)+
R_{\lambda,4}(\phi)+R_{\lambda,5}(\phi)+R_{\lambda,6}(\phi)\bigr)
/ |C_{\lambda}|$.

The set $\fourset{3}(i_{1}, i_{2}, i_{3})$ is a union of six
$M$-sets. Each classical 3-pattern occurs in one of the $M$-sets. 
The cardinality of the $M$-set with pattern $\phi_{c}$ does not
depend on $i_{1}$, $i_{2}$, or $i_{3}$. Hence summing the 
cardinality over all possible $i_{1}$, $i_{2}$, and $i_{3}$ 
is the same as multiplying by $\indfn_{3}(\phi)$. The result is
$R_{\lambda,3}(\phi)$.

The patterns of the $M$-sets in
$\fourset{j}(i_{1}, i_{2}, i_{3})$, $j=4,5,6$, only depend on 
the relative order of $(i_{1}, i_{2}, i_{3}, a_{4}, \ldots, 
a_{j})$. 
Hence, for a fixed such relative order which does not violate $CI_{\phi}$,
the sum of the cardinalities of all $M$-sets with pattern $\phi_{c}$ 
in $\fourset{j}(i_{1}, i_{2}, i_{3})$ is equal to the sum of the 
cardinalities of all $M$-sets with pattern $\phi_{c}$ in 
$\fourset{j}(i_{1}, i_{2}, i_{3})$ 
for fixed $a_{4}, \ldots,a_{j}$ satisfying that relative order, 
multiplied by the number of possible such values of $a_{4}$, \ldots, 
$a_{j}$.

The cardinalities of the $M$-sets do not depend on $i_{1}$, 
$i_{2}$, or $i_{3}$, and summing the number of values of 
$a_{4}, \ldots,a_{j}$ for a given relative order of 
$(i_{1}, i_{2}, i_{3}, a_{4}, \ldots, a_{j})$
over all possible $i_{1}$, $i_{2}$, and $i_{3}$ gives $\indfn_{j}(\phi)$.
Thus $R_{\lambda,j}(\phi)$ can be computed 
by summing the cardinalities of all $M$-sets with pattern $\phi_{c}$ 
in $\fourset{j}(i_{1}, i_{2}, i_{3})$  for
fixed $a_{4}, \ldots, a_{j}$ in each possible relative order of
$(i_{1}, i_{2}, i_{3}, a_{4}, \ldots, a_{j})$, and then multiplying
by $\indfn_{j}(\phi)$.

Suppose that $j \leq n$. Then  
\begin{displaymath}
	\frac{\indfn_{j}(\phi)}{\fall{n}{j}} = 
	\frac{1}{\fall{n}{j}} \binom{n-k_{\phi}}{j-k_{\phi}}  = 
	\frac{\fall{(n-k_{\phi})}{j-k_{\phi}}}{\fall{n}{j}(j-k_{\phi})!} =
	\frac{1}{\fall{n}{k_{\phi}}(j-k_{\phi})!}.
\end{displaymath}

The number $1/\fall{n}{k_{\phi}}(j-k_{\phi})!$
is well defined when $j > n$ too.
The sums $R_{\lambda,j}(\phi)$ are $0$ for all $j > n$, since
the cardinalities of the $M$-sets in $\fourset{j}(i_{1}, i_{2}, i_{3})$
are $0$ when $j > n$. Since $H_{j\ell}(\lambda)=0$ if $j > n$,
all expressions 
$H_{j\ell}(\lambda)/\fall{n}{k_{\phi}}(j-k_{\phi})!=0$ 
when $j > n$.

The above proves the following theorem, which describes
how to compute the means of functions counting vincular
3-patterns as linear combinations of the functions $H_{j\ell}$.

\begin{thm} \label{algo}
	Let $\phi$ be a vincular 3-pattern.
	The mean of $\patt{\phi}$ is given by
	\begin{displaymath}
		\overline{\patt{\phi}} = \frac{1}{\fall{n}{k_{\phi}}}
		\sum_{j=3}^{6} \biggl(\sum_{\ell=1}^{L_{j}}
		u_{j\ell}(\phi) H_{j\ell}(\lambda)\biggr) 
		\frac{1}{(j-k_{\phi})!},
	\end{displaymath}
	where $L_{j}$ is the number of different $M$-set
	cardinalities in $\fourset{j}(i_{1}, i_{2}, i_{3})$, and
	$u_{j\ell}$ is the sum of the number of $M$-sets in 
	$\fourset{j}(i_{1}, i_{2}, i_{3})$ with pattern $\phi_{c}$
	and fixed $a_{4}, \ldots, a_{j}$, for each possible 
	relative order of $(i_{1}, i_{2}, i_{3}, a_{4},\ldots, a_{j})$. 
\end{thm}

To compute $u_{j\ell}(\phi)$, one method is to make a table with
all $M$-sets of cardinality $H_{j\ell} |C_{\lambda}| /
\fall{n}{j}$ on one edge, and every relative order of
$(i_{1}, i_{2}, i_{3}, a_{4},\ldots, a_{j})$ on the other, 
and the patterns of the $M$-sets in the orders in the middle.
Then $u_{j\ell}(\phi)$ is the number of entries $\phi_{c}$ in all
rows with relative orders not violating $CI_{\phi}$.

All $M$-sets $M(a_{4},a_{5},a_{6})$ in 
$\fourset{6}(i_{1}, i_{2}, i_{3})$ have the same cardinality.
Hence $L_{6}=1$ and $u_{61}$ is the number of possible 
relative orders of 
$(i_{1}, i_{2}, i_{3}, a_{4}, a_{5}, a_{6})$ such 
that $M(a_{4},a_{5},a_{6})$ has pattern $\phi_{c}$. 
The pattern $\phi_{c}$ decides the relative order of
$a_{4}$, $a_{5}$, and $a_{6}$. There are four intervals
to put $a_{4}$, $a_{5}$, and $a_{6}$ in: 
$(1,i_{1})$, $(i_{1},i_{2})$, $(i_{2},i_{3})$, and $(i_{3},n)$.
For every condition given by the square brackets and the absence of 
dashes in $\phi$, one of these intervals disappears.
Now $u_{61}$ is equal to the number of ways to put three
identical objects in $4-k_{\phi}$ different boxes (intervals), 
that is $\ttbinom{6-k_{\phi}}{3}$.

Tables for all cardinalities except $H_{61}(\lambda)\frac{|C_{\lambda}|}{
	\fall{n}{6}}$ can be found in the Appendix.

If Theorem~\ref{algo} is used to write $\overline{\patt{\phi}}$
as a linear combination of the functions $H_{j\ell}(\lambda)$, 
then Lemma~\ref{hlincomb} can be used to write $\overline{\patt{\phi}}$
as a linear combination of irreducible characters.

\begin{examp} \label{ex1}
	
	The mean of $\patt{[\text{$21$-$3$})}$ will be computed.
	There are one square bracket and one absent dash in $[\text{$21$-$3$})$. 
	Hence $k_{[\text{$21$-$3$})} = 2$.
	The classical pattern $[\text{$21$-$3$})_{c}$ is $213$, 
	and $CI_{[21-3)}$ consists of the conditions $i_{1}<i_{2}<i_{3}$, \
	$i_{1} = 1$, and $i_{2} = i_{1}+1$.
	
	The relative orders of $(i_{1}, i_{2}, i_{3}, a_{4},
	\ldots, a_{j})$ not
	violating $CI_{[21-3)}$ are exactly those that begin
	with $i_{1}<i_{2}<\ldots$, and they are underlined in the
	tables of the Appendix. Also, all occurrences of the
	pattern $213$ in these rows are underlined.
	
	Now the numbers $u_{j\ell}\bigl([\text{$21$-$3$})\bigr)$ in 
	Theorem~\ref{algo} will be computed.
	
	Table~\ref{mseth11} lists the patterns of all $M$-sets
	of cardinality $H_{31} |C_{\lambda}|/\fall{n}{3}$.
	There are zero occurrences of the pattern $213$ in the
	underlined row in Table~\ref{mseth11}, which means that
	$u_{31}\bigl([\text{$21$-$3$})\bigr) = 0$.
	
	Table~\ref{mseth12} lists the patterns of all $M$-sets
	of cardinality $H_{32} |C_{\lambda}|/\fall{n}{3}$.
	There is one occurrence of the pattern $213$ in the
	underlined row in Table~\ref{mseth12}, which means that
	$u_{32}\bigl([\text{$21$-$3$})\bigr) = 1$.
	
	Similarly, counting the number of occurrences of the 
	pattern $213$ in the underlined rows in 
	Tables~\ref{mseth13}--\ref{mseth32} gives
	\begin{align*}
		u_{33}\bigl([\text{$21$-$3$})\bigr) & = 0, &
		u_{43}\bigl([\text{$21$-$3$})\bigr) & = 2, & 
		u_{52}\bigl([\text{$21$-$3$})\bigr) & = 5. \\
		u_{41}\bigl([\text{$21$-$3$})\bigr) & = 1, &
		u_{44}\bigl([\text{$21$-$3$})\bigr) & = 1, & &\\
		u_{42}\bigl([\text{$21$-$3$})\bigr) & = 2, & 
		u_{51}\bigl([\text{$21$-$3$})\bigr) & = 4, &
	\end{align*}
	
	At last $u_{61}\bigl([\text{$21$-$3$})\bigr)$ is computed as 
	$\ttbinom{6-2}{3} = 4$.
	
	From Theorem~\ref{algo} it now follows that
	\begin{align*}
		\overline{\patt{[\text{$21$-$3$})}} & = \frac{1}{\fall{n}{2}}
		\sum_{j=3}^{6} \biggl(
		\sum_{\ell=1}^{L_{j}} u_{j\ell}\bigl([\text{$21$-$3$})\bigr) 
		H_{j\ell}(\lambda) \biggr)
		\frac{1}{(j-2)!} \\
		& = \frac{1}{n(n-1)} \Bigl(\bigl(0 H_{31}(\lambda)+1 H_{32}(\lambda)+
		0 H_{33}(\lambda) \bigr) \frac{1}{1!} \\
		& \quad + \bigl(1 H_{41}(\lambda)+2 H_{42}(\lambda)+2 H_{43}(\lambda)+
		1 H_{44}(\lambda)\bigr) \frac{1}{2!} \\
		& \quad + \bigl(4 H_{51}(\lambda)+5 H_{52}(\lambda) \bigr) 
		\frac{1}{3!}
		+ 4 H_{61}(\lambda) \frac{1}{4!} \Bigr) \\
		& = \Bigl( H_{32}(\lambda) + \frac{1}{2} H_{41}(\lambda) +
		H_{42}(\lambda)+H_{43}(\lambda) +\frac{1}{2} H_{44}(\lambda) +
		\frac{2}{3} H_{51}(\lambda) \\
		& \quad + \frac{5}{6} H_{52}(\lambda) +\frac{1}{6} H_{61}(\lambda)\Bigr)
		\frac{1}{n(n-1)}.
	\end{align*}
	
	Finally Lemma~\ref{hlincomb} gives that
	
	\begin{align*}
		\overline{\patt{[\text{$21$-$3$})}}& = \biggl(\frac{n-2}{6} 
		\vecsp \frac{(n-3)(n-4)}{6n(n-1)} 
		\vecsp \frac{1}{3n} \vecsp -\frac{n-2}{n(n-1)} \vecsp 0 \vecsp 0 
		\vecsp 0\biggr) \, \chivek.
	\end{align*}
	
\end{examp}

The following theorem gives similar expressions for the rest of
the vincular 3-patterns.

\begin{thm} \label{theotheo}
	The coefficients for the means of all vincular 3-patterns
	expressed as linear combinations of irreducible characters are
	as in Tables~\ref{theo1} and \ref{theo2}, where
	$e_{\ell}$ is the $\ell$:th element in $\chivek$.
	If $n \leq 5$, then some of the elements in $\chivek$ are not
	irreducible characters of $\symg{n}$. The expressions in these cases
	can easily be extracted from Tables~\ref{theo1} and \ref{theo2}
	since $\chi^{(n-2,2)}(\lambda)\!=0$ and $\chi^{(n-3,2,1)}(\lambda) =
	-\chi^{(n-2,1,1)}(\lambda)$ if $n=3$, and $\chi^{(n-3,2,1)}(\lambda)=0$
	if $n=4$.
\end{thm}

{\renewcommand{\textfraction}{0.05}
	\renewcommand{\topfraction}{0.95}
	
	\begin{table}[hptb!]
		{\renewcommand{\arraystretch}{1.5}
			\begin{tabular}{|l|ccccccc|}
				\hline
				Pattern & $e_{1}$ & $e_{2}$ & $e_{3}$ & $e_{4}$
				& $e_{5}$ & $e_{6}$ & $e_{7}$ \\
				\hline
				%
				%
				$(\text{$1$-$2$-$3$})$ & 
				$\frac{1}{6}\ttbinom{n}{3}$ & $\frac{(n+1)(3n-4)}{60}$ & 
				$\frac{n+1}{30}$ & $\frac{3n-2}{30}$ & $0$ & $\frac{1}{15}$ & $0$ \\
				$(\text{$1$-$3$-$2$})$,  
				$(\text{$2$-$1$-$3$})$ &
				$\frac{1}{6}\ttbinom{n}{3}$ & $\frac{(n+1)(n-3)}{60}$ & 
				$-\frac{n+1}{60}$ & $-\frac{n+1}{20}$ & $0$ & $-\frac{1}{30}$ & $0$ \\
				$(\text{$2$-$3$-$1$})$,
				$(\text{$3$-$1$-$2$})$ &
				$\frac{1}{6}\ttbinom{n}{3}$ & $-\frac{(n+1)(3n-4)}{120}$ & 
				$-\frac{n+1}{60}$ & $\frac{n+1}{30}$ & $0$ & $-\frac{1}{30}$ & $0$ \\
				$(\text{$3$-$2$-$1$})$ &
				$\frac{1}{6}\ttbinom{n}{3}$ & $-\frac{(n+1)(n-3)}{30}$ & 
				$\frac{n+1}{30}$ & $-\frac{2n-3}{30}$ & $0$ & $\frac{1}{15}$ & $0$ \\
				\hline
				%
				%
				$(\text{$12$-$3$})$, $(\text{$1$-$23$})$ &
				$\frac{1}{6}\ttbinom{n-1}{2}$ & $\frac{n^2+3n-6}{12n}$ & 
				$\frac{n+2}{12n}$ & $\frac{5n-6}{12n}$ & 0 & $\frac{1}{4n}$ & 0 \\
				$(\text{$13$-$2$})$, $(\text{$2$-$13$})$ &
				$\frac{1}{6}\ttbinom{n-1}{2}$ & $\frac{n-3}{6n}$ & $-\frac{n-1}{6n}$ &  
				$-\frac{1}{2n}$ & 0 & 0 & 0 \\
				$(\text{$21$-$3$})$, $(\text{$1$-$32$})$ &
				$\frac{1}{6}\ttbinom{n-1}{2}$ & $\frac{n(n-3)}{12n}$ & 
				$\frac{n-4}{12n}$ & $-\frac{1}{4}$ & 0 & $-\frac{1}{4n}$ & 0 \\
				$(\text{$23$-$1$})$, $(\text{$3$-$12$})$ &
				$\frac{1}{6}\ttbinom{n-1}{2}$ & $-\frac{n(n-1)}{12n}$ & 
				$\frac{n-4}{12n}$ & $\frac{1}{12}$ & 0 & $-\frac{1}{4n}$ & 0 \\
				$(\text{$31$-$2$})$, $(\text{$2$-$31$})$ &
				$\frac{1}{6}\ttbinom{n-1}{2}$ & $-\frac{2n-3}{6n}$ & 
				$-\frac{n-1}{6n}$ & $\frac{1}{2n}$ & 0 & 0 & 0 \\
				$(\text{$32$-$1$})$, $(\text{$3$-$21$})$ &
				$\frac{1}{6}\ttbinom{n-1}{2}$ & $-\frac{(n+2)(n-3)}{12n}$ & 
				$\frac{n+2}{12n}$ & $-\frac{n-2}{4n}$ & 0 & $\frac{1}{4n}$ & 0 \\
				\hline
				%
				%
				$[\text{$1$-$2$-$3$})$, $(\text{$1$-$2$-$3$}]$  & 
				$\frac{1}{6}\ttbinom{n-1}{2}$ & $\frac{5n-7}{24}$ & $\frac{1}{8}$ & 
				$\frac{5}{24}$ & 0 & $\frac{1}{8n}$ & 0 \\
				$[\text{$1$-$3$-$2$})$, $(\text{$2$-$1$-$3$}]$ & 
				$\frac{1}{6}\ttbinom{n-1}{2}$ & $\frac{n-3}{8}$ & $-\frac{1}{8}$ & 
				$-\frac{5}{24}$ & 0 & $-\frac{1}{8n}$ & 0 \\
				$[\text{$2$-$1$-$3$})$, $(\text{$1$-$3$-$2$}]$  & 
				$\frac{1}{6}\ttbinom{n-1}{2}$ & 0 & 0 & $-\frac{1}{6}$ & 0 & 0 & 0 \\
				$[\text{$2$-$3$-$1$})$, $(\text{$3$-$1$-$2$}]$ & 
				$\frac{1}{6}\ttbinom{n-1}{2}$ & $-\frac{n-2}{6}$ & 0 & 
				$\frac{1}{6}$ & 0 & 0 & 0 \\
				$[\text{$3$-$1$-$2$})$, $(\text{$2$-$3$-$1$}]$ & 
				$\frac{1}{6}\ttbinom{n-1}{2}$ & $-\frac{n+1}{24}$ & $-\frac{1}{8}$ & 
				$\frac{1}{8}$ & 0 & $-\frac{1}{8n}$ & 0 \\
				$[\text{$3$-$2$-$1$})$, $(\text{$3$-$2$-$1$}]$ & 
				$\frac{1}{6}\ttbinom{n-1}{2}$ & $-\frac{n-3}{8}$ & $\frac{1}{8}$ & 
				$-\frac{1}{8}$ & 0 & $\frac{1}{8n}$ & 0 \\
				\hline
				%
				%
				$(\text{$12$-$3$}]$, $[\text{$1$-$23$})$  & 
				$\frac{n-2}{6}$ & $\frac{1}{3}$ & $\frac{n+2}{6\fall{n}{2}}$ & 
				$\frac{1}{2(n-1)}$ & 0 & $\frac{1}{2\fall{n}{2}}$ & 0 \\
				$(\text{$13$-$2$}]$, $[\text{$2$-$13$})$  & 
				$\frac{n-2}{6}$ & $-\frac{n-3}{6n}$ & 
				$-\frac{2n-5}{6\fall{n}{2}}$ & $-\frac{1}{2\fall{n}{2}}$ & 0 & 
				$\frac{1}{2\fall{n}{2}}$ & 0 \\
				$(\text{$21$-$3$}]$, $[\text{$1$-$32$})$  & 
				$\frac{n-2}{6}$ & $\frac{(2n+1)(n-3)}{6\fall{n}{2}}$ & 
				$\frac{n-7}{6\fall{n}{2}}$ & $-\frac{n+1}{2\fall{n}{2}}$ & 
				0 & $-\frac{1}{\fall{n}{2}}$ & 0 \\
				$(\text{$23$-$1$}]$, $[\text{$3$-$12$})$  & 
				$\frac{n-2}{6}$ & $-\frac{n^2-4n+9}{6\fall{n}{2}}$ & 
				$\frac{n-7}{6\fall{n}{2}}$ & $\frac{n-3}{2\fall{n}{2}}$ & 
				0 & $-\frac{1}{\fall{n}{2}}$ & 0 \\
				$(\text{$31$-$2$}]$, $[\text{$2$-$31$})$  & 
				$\frac{n-2}{6}$ & $-\frac{n^2+2n-9}{6\fall{n}{2}}$ & 
				$-\frac{2n-5}{6\fall{n}{2}}$ & $\frac{3}{2\fall{n}{2}}$ & 
				0 & $\frac{1}{2\fall{n}{2}}$ & 0 \\
				$(\text{$32$-$1$}]$, $[\text{$3$-$21$})$ &
				$\frac{n-2}{6}$ & $-\frac{(n+2)(n-3)}{6\fall{n}{2}}$ & 
				$\frac{n+2}{6\fall{n}{2}}$ & $-\frac{n-2}{2\fall{n}{2}}$ & 
				0 & $\frac{1}{2\fall{n}{2}}$ & 0 \\
				\hline
				%
				%
				$[\text{$12$-$3$})$, $(\text{$1$-$23$}]$ & 
				$\frac{n-2}{6}$ & $\frac{n^2+5n-12}{6\fall{n}{2}}$ & 
				$\frac{1}{3n}$ & $\frac{n-2}{\fall{n}{2}}$ & 0 & 0 & 0 \\
				$[\text{$13$-$2$})$, $(\text{$2$-$13$}]$  & 
				$\frac{n-2}{6}$ & $\frac{(n+2)(n-3)}{6\fall{n}{2}}$ & 
				$-\frac{2}{3n}$ & $-\frac{1}{\fall{n}{2}}$ & 0 & 0 & 0 \\
				$[\text{$21$-$3$})$, $(\text{$1$-$32$}]$  & 
				$\frac{n-2}{6}$ & $\frac{(n-3)(n-4)}{6\fall{n}{2}}$ & 
				$\frac{1}{3n}$ & $-\frac{n-2}{\fall{n}{2}}$ & 0 & 0 & 0 \\
				$[\text{$23$-$1$})$, $(\text{$3$-$12$}]$  & 
				$\frac{n-2}{6}$ & $-\frac{n^2-n-3}{3\fall{n}{2}}$ & 
				$\frac{1}{3n}$ & $\frac{1}{\fall{n}{2}}$ & 0 & 0 & 0 \\
				$[\text{$31$-$2$})$, $(\text{$2$-$31$}]$  & 
				$\frac{n-2}{6}$ & $\frac{n-6}{6n}$ & $-\frac{2}{3n}$ & 
				$\frac{1}{\fall{n}{2}}$ & 0 & 0 & 0 \\
				$[\text{$32$-$1$})$, $(\text{$3$-$21$}]$, $[\text{$3$-$2$-$1$}]$  & 
				$\frac{n-2}{6}$ & $-\frac{n-3}{3n}$ & $\frac{1}{3n}$ & 
				$-\frac{1}{\fall{n}{2}}$ & 0 & 0 & 0 \\
				\hline
				%
				%
				$[\text{$1$-$2$-$3$}]$  & 
				$\frac{n-2}{6}$ & $\frac{2n-3}{3n}$ & $\frac{1}{3n}$ & 
				$\frac{1}{\fall{n}{2}}$ & 0 & 0 & 0 \\
				$[\text{$1$-$3$-$2$}]$, $[\text{$2$-$1$-$3$}]$  & 
				$\frac{n-2}{6}$ & $\frac{n-3}{6n}$ & $-\frac{1}{6n}$ & $-\frac{1}{2n}$ & 
				0 & 0 & 0 \\
				$[\text{$2$-$3$-$1$}]$, $[\text{$3$-$1$-$2$}]$ &
				$\frac{n-2}{6}$ & $-\frac{2n-3}{6n}$ & $-\frac{1}{6n}$ & $\frac{1}{2n}$ & 
				0 & 0 & 0 \\
				\hline
		\end{tabular}}
		\caption{Coefficients for means of 
			vincular 3-patterns.\newline Note that 
			$\overline{\patt{(\text{$32$-$1$})}} = 
			\frac{n-1}{2}\overline{\patt{(\text{$32$-$1$}]}}$.}
		\label{theo1}
	\end{table}
	
	\begin{table}[ht!]
		{\renewcommand{\arraystretch}{1.5}
			\begin{tabular}{|l|ccccccc|}
				\hline
				Pattern & $e_{1}$ & $e_{2}$ & $e_{3}$ & $e_{4}$
				& $e_{5}$ & $e_{6}$ & $e_{7}$ \\
				\hline
				%
				%
				$(\text{$123$})$ &
				$\frac{n-2}{6}$ & $\frac{1}{n}$ & $\frac{1}{\fall{n}{2}}$ & 
				$\frac{1}{n}$ & 0 & $\frac{1}{\fall{n}{2}}$ & 0 \\
				$(\text{$132$})$, $(\text{$213$})$ &
				$\frac{n-2}{6}$ & $\frac{n-3}{2\fall{n}{2}}$ & $-\frac{1}{2\fall{n}{2}}$ & 
				$-\frac{3}{2\fall{n}{2}}$ & 0 & $-\frac{1}{2\fall{n}{2}}$ & 0 \\
				$(\text{$231$})$, $(\text{$312$})$ &
				$\frac{n-2}{6}$ & $-\frac{1}{2n}$ & $-\frac{1}{2\fall{n}{2}}$ & 
				$\frac{1}{2\fall{n}{2}}$ & 0 & $-\frac{1}{2\fall{n}{2}}$ & 0 \\
				$(\text{$321$})$ &
				$\frac{n-2}{6}$ & $-\frac{n-3}{\fall{n}{2}}$ & $\frac{1}{\fall{n}{2}}$ & 
				$-\frac{n-3}{\fall{n}{2}}$ & 0 & $\frac{1}{\fall{n}{2}}$ & 0 \\
				\hline
				%
				%
				$[\text{$123$})$, $(\text{$123$}]$  &
				$\frac{1}{6}$ & $\frac{1}{n(n-2)}$ & $\frac{1}{\fall{n}{3}}$ & 
				$\frac{1}{n(n-2)}$ & 0 & $\frac{1}{\fall{n}{3}}$ & 0 \\
				$[\text{$132$})$, $(\text{$213$}]$ &
				$\frac{1}{6}$ & $\frac{2n-6}{\fall{n}{3}}$ & $-\frac{2}{\fall{n}{3}}$ & 
				$\frac{n-6}{\fall{n}{3}}$ & 0 & $-\frac{2}{\fall{n}{3}}$ & 0 \\
				$[\text{$213$})$, $(\text{$132$}]$, $[\text{$321$})$, $(\text{$321$}]$ &
				$\frac{1}{6}$ & $-\frac{n-3}{\fall{n}{3}}$ & $\frac{1}{\fall{n}{3}}$ & 
				$-\frac{n-3}{\fall{n}{3}}$ & 0 & $\frac{1}{\fall{n}{3}}$ & 0 \\
				$[\text{$231$})$, $(\text{$312$}]$ &
				$\frac{1}{6}$ & $-\frac{2n-5}{\fall{n}{3}}$ & $\frac{1}{\fall{n}{3}}$ & 
				$-\frac{n-5}{\fall{n}{3}}$ & 0 & $\frac{1}{\fall{n}{3}}$ & 0 \\
				$[\text{$312$})$, $(\text{$231$}]$ &
				$\frac{1}{6}$ & $\frac{n-4}{\fall{n}{3}}$ & $-\frac{2}{\fall{n}{3}}$ & 
				$\frac{n-4}{\fall{n}{3}}$ & 0 & $-\frac{2}{\fall{n}{3}}$ & 0 \\
				\hline
				%
				%
				$[\text{$12$-$3$}]$, $[\text{$1$-$23$}]$  &
				$\frac{1}{6}$ & $\frac{n+1}{2\fall{n}{2}}$ & $\frac{1}{2\fall{n}{2}}$ & 
				$\frac{1}{2\fall{n}{2}}$ & 0 & $-\frac{1}{2\fall{n}{3}}$ & 0 \\
				$[\text{$13$-$2$}]$, $[\text{$2$-$13$}]$ & 
				$\frac{1}{6}$ & 0 & $-\frac{1}{\fall{n}{2}}$ & 0 & 0 & 
				$\frac{1}{\fall{n}{3}}$ & 0 \\
				$[\text{$21$-$3$}]$, $[\text{$1$-$32$}]$ &
				$\frac{1}{6}$ & $\frac{n-3}{2\fall{n}{2}}$ & $\frac{1}{2\fall{n}{2}}$ & 
				$-\frac{3}{2\fall{n}{2}}$ & 0 & $-\frac{1}{2\fall{n}{3}}$ & 0 \\
				$[\text{$23$-$1$}]$, $[\text{$3$-$12$}]$ & 
				$\frac{1}{6}$ & $-\frac{1}{2n}$ & $\frac{1}{2\fall{n}{2}}$ & 
				$\frac{1}{2\fall{n}{2}}$ & 0 & $-\frac{1}{2\fall{n}{3}}$ & 0 \\
				$[\text{$31$-$2$}]$, $[\text{$2$-$31$}]$ & 
				$\frac{1}{6}$ & $-\frac{1}{\fall{n}{2}}$ & $-\frac{1}{\fall{n}{2}}$ & 
				$\frac{2}{\fall{n}{3}}$ & 0 & $\frac{1}{\fall{n}{3}}$ & 0 \\
				$[\text{$32$-$1$}]$, $[\text{$3$-$21$}]$ & 
				$\frac{1}{6}$ & $-\frac{n-3}{2\fall{n}{2}}$ & $\frac{1}{2\fall{n}{2}}$ & 
				$\frac{n-6}{2\fall{n}{3}}$ & 0 & $-\frac{1}{2\fall{n}{3}}$ & 0 \\
				\hline
		\end{tabular}}
		\caption{Coefficients for means of 
			vincular 3-patterns. Note that \newline 
			$\overline{\patt{(\text{$123$})}} = 
			(n-2)\overline{\patt{[\text{$123$})}}$ and that
			$\overline{\patt{(\text{$321$})}} = 
			(n-2)\overline{\patt{[\text{$321$})}}$.}
		\label{theo2}
	\end{table}

	\begin{bevis}
		The results in the tables follows from Theorem~\ref{algo} and
		Tables~\ref{mseth11}--\ref{mseth32} of the Appendix.
	\end{bevis}
	
	As seen in Tables~\ref{theo1} and \ref{theo2}, there are
	many pairs of means with identical coefficients. This is due to 
	a bijection $\psi:\symg{n} \to \symg{n}$, described in \cite{simschm},
	whose restriction to each conjugacy class $C_{\lambda}$ is a
	bijection $\psi_{|C_{\lambda}}:C_{\lambda} \to C_{\lambda}$.
	Let $\sigma \in S_{n}$ be the permutation $n \; (n-1)\;\ldots\;2\;1$.
	The bijection $\psi$ is defined by
	$\psi(\pi) = \sigma^{-1} \pi \sigma$. 
	If $\tau = \psi(\pi)$, then $\tau_{i} = n+1-\pi_{n+1-i}$. Consider the
	sequences $\pi_{i_{1}} \pi_{i_{2}} \pi_{i_{3}}$ and
	$\tau_{n+1-i_{3}} \tau_{n+1-i_{2}} \tau_{n+1-i_{1}}$, where $i_{1}<i_{2}<i_{3}$. 
	It is not hard to deduce that the first sequence is an occurrence
	of a pattern $\phi$ in Table~\ref{theo1} or \ref{theo2}
	if and only if the second
	sequence is an occurrence of its ``pair'' $\phi'$. The pattern
	$\phi'$ can be constructed from $\phi$ by reversing the positions
	of all brackets and dashes, and by the following equivalences
	between $\phi_{c}$ and $\phi'_{c}$:
	\begin{align*}
		123 & \leftrightarrow 123, &
		132 & \leftrightarrow 213, &
		231 & \leftrightarrow 312, &
		321 & \leftrightarrow 321.
	\end{align*}

	Many well-studied permutation statistics are linear combinations
	of pattern-counting statistics. Hence they, too, can be analysed
	using Theorem~\ref{theotheo}. One example is the number of peaks, 
	which is mentioned as a corollary.
	
	\begin{cor} \label{peakth}
		The mean $\overline{\peak}$ can be written
		\begin{displaymath}
			\overline{\peak} = \frac{n-2}{3} \chi^{(n)}-\frac{1}{n(n-1)}
			\bigl( \chi^{(n-1,1)}+ \chi^{(n-2,2)}+ \chi^{(n-2,1,1)}+ 
			\chi^{(n-3,2,1)} \bigr)
		\end{displaymath}
	\end{cor}
	
	\begin{bevis}
		Follows from the fact that $\overline{\peak} = 
		\overline{\patt{(132)}} + \overline{\patt{(231)}}$ and 
		Theorem~\ref{theotheo}.
	\end{bevis}
	
	The procedure described in this section can be generalised to
	patterns of any length. It is obvious that the cardinalities of the 
	generalised $M$-sets will still be class functions. Hence
	all $H_{j\ell}(\lambda)$ will be class functions and linear
	combinations of the irreducible characters in $\symg{n}$,
	since these characters form a basis for all class functions. It seems
	like all irreducible characters needed in the case of $k$-patterns
	are those that are indexed by partitions 
	$\lambda = (\lambda_1, \lambda_2,\ldots, \lambda_{\ell})$ where
	$\lambda_1 \geq (n-k)$.
	The arguments given in this section also apply with minor 
	changes to the general case. The only drawback is that the 
	number of functions $H_{j\ell}(\lambda)$ and tables for
	computing the numbers $u_{j\ell}$ will 
	increase, and so will the size of these tables.
}

\section{Expected values} \label{expvals}

As seen in Theorem~\ref{theotheo}, if $\phi$ is a (classical or)
vincular 3-pattern, then the mean
$\overline{\patt{\phi}}$
is a linear combination of five of the seven irreducible characters
in Lemma~\ref{chin}.

Let 
$ \subsetpn =  \{(n), (n-1,1), (n-2,2), (n-2,1,1), (n-3,2,1)\} 
\subset P_{n}$.

It is now easy
to compute the expected value $\mathbb{E}_{\Gamma}(\patt{\phi},t)$
as described in Section~\ref{theory}. If 
$\overline{\patt{\phi}} = \sum_{\lambda \vdash n} a_{\lambda} \chi^{\lambda}$, 
then $a_{\lambda}=0$ for all $\lambda \notin \subsetpn$.
Together with Theorem~\ref{hultone}, this gives
\begin{displaymath}
	\mathbb{E}_{\Gamma}(\patt{\phi},t) = \frac{n!}{|\Gamma|^t} \sum_{\lambda \in 
		\subsetpn} 
	a_{\lambda} b_{\lambda}^{(t)},
\end{displaymath}
where the coefficients $b_{\lambda}^{(t)}$ are given by 
Theorem~\ref{hulttwo} when $\Gamma$ is a conjugacy class.
By Lemma~\ref{chin} it is very easy to calculate the five
needed coefficients $b_{\lambda}^{(t)}$ explicitly.
If $\Gamma = C_{\mu}$ and if $\mu$ has $p$ parts of size $1$, 
$q$ parts of size $2$, and $r$ parts of size $3$, then 
\begin{align*}
	b_{(n)}^{(t)} &= 
	\frac{|C_{\mu}|^t}{n!},\\
	b_{(n-1,1)}^{(t)} &= \frac{|C_{\mu}|^t}{n!} \frac{(p-1)^t}{(n-1)^{t-1}},\\
	b_{(n-2,2)}^{(t)} &= \frac{|C_{\mu}|^t}{n!} \frac{\bigl(p(p-3
		)/2+q\bigr)^t}{\bigl(n(n-3)/2\bigr)^{t-1}},\\
	b_{(n-2,1,1)}^{(t)} &= \frac{|C_{\mu}|^t}{n!} \frac{\bigl((p-1)(p-2)/2-q
		\bigr)^t}{\bigl((n-1)(n-2)/2\bigr)^{t-1}},\\
	b_{(n-3,2,1)}^{(t)} &= \frac{|C_{\mu}|^t}{n!} \frac{\bigl(p(p-2)(p-4)/3
		-r\bigr)^t}{\bigl(n(n-2)(n-4)/3\bigr)^{t-1}}. 
\end{align*}

{\renewcommand{\textfraction}{0.05}
	\renewcommand{\topfraction}{0.95}
	\renewcommand{\bottomfraction}{0.95}
	\setcounter{totalnumber}{5}
	\setcounter{topnumber}{5}
	\setcounter{bottomnumber}{5}
	
	\begin{examp} Suppose that $\Gamma = T$, the set of transpositions
		in $\symg{n}$. The cycle type of the transpositions has $n-2$ 
		parts of size $1$ and one part of size $2$, and $|T| = \ttbinom{n}{2}
		= \frac{n(n-1)}{2}$.
		Hence the numbers $b_{\lambda}^{(t)}$ are given by
		\begin{align*}
			b_{(n)}^{(t)} &= \frac{1}{n!}\frac{(n-1)^t n^t}{2^t},\\
			b_{(n-1,1)}^{(t)} &=  \frac{1}{n!} \frac{(n-3)^t n^t}{2^t}(n-1),\\
			b_{(n-2,2)}^{(t)} &= 
			\frac{1}{n!}\frac{(n-1)^t (n-4)^t}{2^t} \frac{n(n-3)}{2},\\
			b_{(n-2,1,1)}^{(t)} &= \frac{1}{n!} \frac{(n-5)^t n^t}{2^t}\frac{(n-1)(n-2)}{2},\\
			b_{(n-3,2,1)}^{(t)} &= \frac{1}{n!} \frac{(n-1)^t (n-6)^t}{2^t}
			\frac{n(n-2)(n-4)}{3}.
		\end{align*}

		The expected values of the number of occurrences of the classical
		3-patterns after $t$ random transpositions can now be computed
		using the expressions in Theorem~\ref{theotheo}. They are 
		\begin{align*}
			\begin{split}
				& \mathbb{E}_{T}(\patt{(\text{$1$-$2$-$3$})},t) = \frac{1}{6}\binom{n}{3}+
				\frac{(n+1)(n-1)(3n-4)}{60}\left(1-\frac{2}{n-1}\right)^t \\
				& \qquad + \frac{(n+1)n(n-3)}{60}\left(1-\frac{4}{n}\right)^t +
				\frac{(3n-2)(n-1)(n-2)}{60}\left(1-\frac{4}{n-1}\right)^t \\ 
				& \qquad + \frac{n(n-2)(n-4)}{45}\left(1-\frac{6}{n}\right)^t,
			\end{split} \\
			\begin{split}
				& \mathbb{E}_{T}(\patt{(\text{$3$-$2$-$1$})},t) = \frac{1}{6}\binom{n}{3}-
				\frac{(n+1)(n-1)(n-3)}{30}\left(1-\frac{2}{n-1}\right)^t \\
				& \qquad + \frac{(n+1)n(n-3)}{60}\left(1-\frac{4}{n}\right)^t +
				\frac{(2n-3)(n-1)(n-2)}{60}\left(1-\frac{4}{n-1}\right)^t \\ 
				& \qquad + \frac{n(n-2)(n-4)}{45}\left(1-\frac{6}{n}\right)^t,
			\end{split} \\
			\begin{split}
				& \mathbb{E}_{T}(\patt{(\text{$1$-$3$-$2$})},t) = 
				\mathbb{E}_{T}(\patt{(\text{$2$-$1$-$3$})},t) = \frac{1}{6}\binom{n}{3}+
				\frac{(n+1)(n-1)(n-3)}{60}\left(1-\frac{2}{n-1}\right)^t \\
				& \qquad - \frac{(n+1)n(n-3)}{120}\left(1-\frac{4}{n}\right)^t -
				\frac{(n+1)(n-1)(n-2)}{40}\left(1-\frac{4}{n-1}\right)^t \\ 
				& \qquad - \frac{n(n-2)(n-4)}{90}\left(1-\frac{6}{n}\right)^t,
			\end{split} \\
			\begin{split}
				& \mathbb{E}_{T}(\patt{(\text{$2$-$3$-$1$})},t) = 
				\mathbb{E}_{T}(\patt{(\text{$3$-$1$-$2$})},t) = \frac{1}{6}\binom{n}{3}-
				\frac{(n+1)(n-1)(3n-4)}{120}\left(1-\frac{2}{n-1}\right)^t \\
				& \qquad - \frac{(n+1)n(n-3)}{120}\left(1-\frac{4}{n}\right)^t +
				\frac{(n+1)(n-1)(n-2)}{60}\left(1-\frac{4}{n-1}\right)^t \\ 
				& \qquad - \frac{n(n-2)(n-4)}{90}\left(1-\frac{6}{n}\right)^t.
			\end{split}
		\end{align*}
		
	\end{examp}

	{\setcounter{table}{0}
		\renewcommand{\thetable}{A.\arabic{table}}
		
		\addcontentsline{toc}{section}{Appendix}
		\section*{Appendix\markright{Appendix}{}} \label{app}

		\begin{table}[htp!]
			\begin{tabular}{|c|c|}
				\hline
				& \Mset{i_{1}}{i_{2}}{i_{3}}  \\
				\hline
				$\underline{i_{1}<i_{2}<i_{3}}$ & 123 \\[0.07cm]
				\hline
			\end{tabular} 
			\caption{$M$-sets of cardinality $H_{31}(\lambda)
				|C_{\lambda}|/\fall{n}{3}$. \newline
				Underlines refer to Example~\ref{ex1}.}
			\label{mseth11}
		\end{table}
		
		\begin{table}[htp!]
			\begin{tabular}{|c|c|c|c|}
				\hline
				& \Mset{i_{1}}{i_{3}}{i_{2}} & \Mset{i_{3}}{i_{2}}{i_{1}} & 
				\Mset{i_{2}}{i_{1}}{i_{3}} \\
				\hline
				$\underline{i_{1}<i_{2}<i_{3}}$ & 132 & 321 & \underline{213} \\[0.07cm]
				\hline
			\end{tabular} 
			\caption{$M$-sets of cardinality $H_{32}(\lambda)
				|C_{\lambda}|/\fall{n}{3}$. \newline
				Underlines refer to Example~\ref{ex1}.}
			\label{mseth12}
		\end{table}
		
		\begin{table}[htp!]
			\begin{tabular}{|c|c|c|}
				\hline
				& \Mset{i_{2}}{i_{3}}{i_{1}} & \Mset{i_{3}}{i_{1}}{i_{2}}  \\
				\hline
				$\underline{i_{1}<i_{2}<i_{3}}$ & 231 & 312 \\[0.07cm]
				\hline
			\end{tabular} 
			\caption{$M$-sets of cardinality $H_{33}(\lambda)
				|C_{\lambda}|/\fall{n}{3}$. \newline
				Underlines refer to Example~\ref{ex1}.}
			\label{mseth13}
		\end{table}

		\begin{table}[htp!]
			\begin{tabular}{|c|c|c|c|}
				\hline
				& \Mset{i_{1}}{i_{2}}{a_{4}} & \Mset{i_{1}}{a_{4}}{i_{3}} & \Mset{a_{4}}{i_{2}}{i_{3}}  \\
				\hline
				$a_{4}<i_{1}<i_{2}<i_{3}$ & 231 & 213 & 123 \\
				$i_{1}<a_{4}<i_{2}<i_{3}$ & 132 & 123 & 123 \\
				$\underline{i_{1}<i_{2}<a_{4}<i_{3}}$ & 123 & 123 & \underline{213} \\
				$\underline{i_{1}<i_{2}<i_{3}<a_{4}}$ & 123 & 132 & 312 \\[0.07cm]
				\hline
			\end{tabular} 
			\caption{$M$-sets of cardinality $H_{41}(\lambda)
				|C_{\lambda}|/\fall{n}{4}$. \newline
				Underlines refer to Example~\ref{ex1}.}
			\label{mseth21}
		\end{table}
		
		\begin{table}[htp!]
			\begin{tabular}{|c|c|c|c|}
				\hline
				& \Mset{i_{2}}{i_{1}}{a_{4}} & \Mset{i_{3}}{a_{4}}{i_{1}} & \Mset{a_{4}}{i_{3}}{i_{2}}\\
				\hline
				$a_{4}<i_{1}<i_{2}<i_{3}$ & 321 & 312 & 132 \\
				$i_{1}<a_{4}<i_{2}<i_{3}$ & 312 & 321 & 132 \\
				$\underline{i_{1}<i_{2}<a_{4}<i_{3}}$ & \underline{213} & 321 & 231 \\
				$\underline{i_{1}<i_{2}<i_{3}<a_{4}}$ & \underline{213} & 231 & 321 \\[0.07cm]
				\hline
			\end{tabular} 
			\caption{$M$-sets of cardinality $H_{42}(\lambda)
				|C_{\lambda}|/\fall{n}{4}$. \newline
				Underlines refer to Example~\ref{ex1}.}
			\label{mseth22}
		\end{table}

		\begin{table}[htp!]
			\begin{tabular}{|c|c|c|c|}
				\hline
				& \Mset{i_{1}}{i_{3}}{a_{4}} & \Mset{i_{1}}{a_{4}}{i_{2}} & \Mset{i_{3}}{i_{2}}{a_{4}} \\
				\hline
				$a_{4}<i_{1}<i_{2}<i_{3}$ & 231 & 213 & 321 \\
				$i_{1}<a_{4}<i_{2}<i_{3}$ & 132 & 123 & 321 \\
				$\underline{i_{1}<i_{2}<a_{4}<i_{3}}$ & 132 & 132 & 312  \\
				$\underline{i_{1}<i_{2}<i_{3}<a_{4}}$ & 123 & 132 & \underline{213} \\[0.07cm]
				\hline
				&
				\Mset{a_{4}}{i_{2}}{i_{1}} & \Mset{i_{2}}{a_{4}}{i_{3}} & \Mset{a_{4}}{i_{1}}{i_{3}}\\
				\hline
				$a_{4}<i_{1}<i_{2}<i_{3}$  & 132 & 213 & 123 \\
				$i_{1}<a_{4}<i_{2}<i_{3}$  & 231 & 213 & 213 \\
				$\underline{i_{1}<i_{2}<a_{4}<i_{3}}$ & 321 & 123 & \underline{213} \\
				$\underline{i_{1}<i_{2}<i_{3}<a_{4}}$ & 321 & 132 & 312 \\[0.07cm]
				\hline
			\end{tabular}
			\caption{$M$-sets of cardinality $H_{43}(\lambda)
				|C_{\lambda}|/\fall{n}{4}$. \newline
				Underlines refer to Example~\ref{ex1}.}
			\label{mseth23}
		\end{table}

		\begin{table}[htp!]
			\begin{tabular}{|c|c|c|c|}
				\hline
				& \Mset{i_{2}}{i_{3}}{a_{4}} & \Mset{i_{3}}{a_{4}}{i_{2}} & \Mset{i_{3}}{i_{1}}{a_{4}} \\
				\hline
				$a_{4}<i_{1}<i_{2}<i_{3}$ & 231 & 312 & 321 \\
				$i_{1}<a_{4}<i_{2}<i_{3}$ & 231 & 312 & 312  \\
				$\underline{i_{1}<i_{2}<a_{4}<i_{3}}$ & 132 & 321 & 312  \\
				$\underline{i_{1}<i_{2}<i_{3}<a_{4}}$ & 123 & 231 & \underline{213} \\[0.07cm]
				\hline
				&
				\Mset{a_{4}}{i_{3}}{i_{1}} & \Mset{i_{2}}{a_{4}}{i_{1}} & \Mset{a_{4}}{i_{1}}{i_{2}}\\
				\hline
				$a_{4}<i_{1}<i_{2}<i_{3}$ & 132 & 312 & 123 \\
				$i_{1}<a_{4}<i_{2}<i_{3}$ & 231 & 321 & 213 \\
				$\underline{i_{1}<i_{2}<a_{4}<i_{3}}$ & 231 & 231 & 312 \\
				$\underline{i_{1}<i_{2}<i_{3}<a_{4}}$ & 321 & 231 & 312 \\[0.07cm]
				\hline
			\end{tabular}
			\caption{$M$-sets of cardinality $H_{44}(\lambda)
				|C_{\lambda}|/\fall{n}{4}$. \newline
				Underlines refer to Example~\ref{ex1}.}
			\label{mseth24}
		\end{table}
		
		\newpage
		Note that in the last three pairs of rows of Tables~\ref{mseth31} and 
		\ref{mseth32} all six classical 3-patterns occur exactly the
		same number of times. The reason is that all of $i_{1}$, $i_{2}$, and 
		$i_{3}$ have the same relation to $a_{4}$ and $a_{5}$
		in these relative orders, that the relative order of $a_{4}$ and 
		$a_{5}$ is shifted
		in each pair of rows, and that the indices $i_{\ell}$ in the $M$-sets have
		three different positions once (in Table~\ref{mseth31}) or twice
		(in Table~\ref{mseth32}). Hence these rows could be
		omitted and some constant, depending on the number
		of square brackets and dashes in $\phi$, could be added to each $u_{j\ell}$, 
		but for simplicity they are included.

		\begin{table}[htp!]
			\begin{tabular}{|l|c|c|c|}
				\hline
				& \Mset{i_{1}}{a_{4}}{a_{5}} & \Mset{a_{4}}{i_{2}}{a_{5}} & \Mset{a_{4}}{a_{5}}{i_{3}}\\
				\hline
				$a_{4}<i_{1}<a_{5}<i_{2}<i_{3}$ & 213 & 132 & 123 \\
				$a_{5}<i_{1}<a_{4}<i_{2}<i_{3}$ & 231 & 231 & 213 \\
				$a_{4}<i_{1}<i_{2}<a_{5}<i_{3}$ & 213 & 123 & 123 \\
				$a_{5}<i_{1}<i_{2}<a_{4}<i_{3}$ & 231 & 321 & 213 \\
				$i_{1}<a_{4}<a_{5}<i_{2}<i_{3}$ & 123 & 132 & 123 \\
				$i_{1}<a_{5}<a_{4}<i_{2}<i_{3}$ & 132 & 231 & 213 \\
				$i_{1}<a_{4}<i_{2}<a_{5}<i_{3}$ & 123 & 123 & 123 \\
				$i_{1}<a_{5}<i_{2}<a_{4}<i_{3}$ & 132 & 321 & 213 \\
				$i_{1}<a_{4}<i_{2}<i_{3}<a_{5}$ & 123 & 123 & 132 \\
				$i_{1}<a_{5}<i_{2}<i_{3}<a_{4}$ & 132 & 321 & 312 \\
				$\underline{i_{1}<i_{2}<a_{4}<a_{5}<i_{3}}$ & 123 & \underline{213} & 123 \\
				$\underline{i_{1}<i_{2}<a_{5}<a_{4}<i_{3}}$ & 132 & 312 & \underline{213} \\
				$\underline{i_{1}<i_{2}<a_{4}<i_{3}<a_{5}}$ & 123 & \underline{213} & 132 \\
				$\underline{i_{1}<i_{2}<a_{5}<i_{3}<a_{4}}$ & 132 & 312 & 312 \\[0.07cm]
				\hline
				$a_{4}<a_{5}<i_{1}<i_{2}<i_{3}$ & 312 & 132 & 123 \\
				$a_{5}<a_{4}<i_{1}<i_{2}<i_{3}$ & 321 & 231 & 213 \\
				$a_{4}<i_{1}<i_{2}<i_{3}<a_{5}$ & 213 & 123 & 132 \\
				$a_{5}<i_{1}<i_{2}<i_{3}<a_{4}$ & 231 & 321 & 312 \\
				$\underline{i_{1}<i_{2}<i_{3}<a_{4}<a_{5}}$ & 123 & \underline{213} & 231 \\
				$\underline{i_{1}<i_{2}<i_{3}<a_{5}<a_{4}}$ & 132 & 312 & 321 \\[0.07cm]
				\hline
			\end{tabular}
			\caption{$M$-sets of cardinality $H_{51}(\lambda)
				|C_{\lambda}|/\fall{n}{5}$. \newline
				Underlines refer to Example~\ref{ex1}.}
			\label{mseth31}
		\end{table}

		\begin{table}[htp!]
			\begin{tabular}{|c|c|c|c|}
				\hline
				& \Mset{i_{2}}{a_{4}}{a_{5}} & \Mset{i_{3}}{a_{4}}{a_{5}} & \Mset{a_{4}}{i_{1}}{a_{5}} \\
				\hline
				$ a_{4}<i_{1}<a_{5}<i_{2}<i_{3}$ & 312 & 312 & 123  \\
				$ a_{5}<i_{1}<a_{4}<i_{2}<i_{3}$ & 321 & 321 & 321  \\
				$ a_{4}<i_{1}<i_{2}<a_{5}<i_{3}$ & 213 & 312 & 123  \\
				$ a_{5}<i_{1}<i_{2}<a_{4}<i_{3}$ & 231 & 321 & 321  \\
				$ i_{1}<a_{4}<a_{5}<i_{2}<i_{3}$ & 312 & 312 & 213  \\
				$ i_{1}<a_{5}<a_{4}<i_{2}<i_{3}$ & 321 & 321 & 312  \\
				$ i_{1}<a_{4}<i_{2}<a_{5}<i_{3}$ & 213 & 312 & 213  \\
				$ i_{1}<a_{5}<i_{2}<a_{4}<i_{3}$ & 231 & 321 & 312  \\
				$ i_{1}<a_{4}<i_{2}<i_{3}<a_{5}$ & 213 & 213 & 213  \\
				$ i_{1}<a_{5}<i_{2}<i_{3}<a_{4}$ & 231 & 231 & 312  \\
				$\underline{ i_{1}<i_{2}<a_{4}<a_{5}<i_{3}}$ & 123 & 312 & \underline{213} \\
				$\underline{ i_{1}<i_{2}<a_{5}<a_{4}<i_{3}}$ & 132 & 321 & 312 \\
				$\underline{ i_{1}<i_{2}<a_{4}<i_{3}<a_{5}}$ & 123 & \underline{213} & 
				\underline{213} \\
				$\underline{ i_{1}<i_{2}<a_{5}<i_{3}<a_{4}}$ & 132 & 231 & 312 \\[0.07cm]
				\hline
				$a_{4}<a_{5}<i_{1}<i_{2}<i_{3}$ & 312 & 312 & 132  \\
				$a_{5}<a_{4}<i_{1}<i_{2}<i_{3}$ & 321 & 321 & 231  \\
				$a_{4}<i_{1}<i_{2}<i_{3}<a_{5}$ & 213 & 213 & 123  \\
				$a_{5}<i_{1}<i_{2}<i_{3}<a_{4}$ & 231 & 231 & 321  \\
				$\underline{i_{1}<i_{2}<i_{3}<a_{4}<a_{5}}$ & 123 & 123 & \underline{213}  \\
				$\underline{i_{1}<i_{2}<i_{3}<a_{5}<a_{4}}$ & 132 & 132 & 312 \\[0.07cm]
				\hline
				&
				\Mset{a_{4}}{i_{3}}{a_{5}} & \Mset{a_{4}}{a_{5}}{i_{1}} & \Mset{a_{4}}{a_{5}}{i_{2}}\\
				\hline
				$ a_{4}<i_{1}<a_{5}<i_{2}<i_{3}$ & 132 & 132 & 123 \\
				$ a_{5}<i_{1}<a_{4}<i_{2}<i_{3}$ & 231 & 312 & 213 \\
				$ a_{4}<i_{1}<i_{2}<a_{5}<i_{3}$ & 132 & 132 & 132 \\
				$ a_{5}<i_{1}<i_{2}<a_{4}<i_{3}$ & 231 & 312 & 312 \\
				$ i_{1}<a_{4}<a_{5}<i_{2}<i_{3}$ & 132 & 231 & 123 \\
				$ i_{1}<a_{5}<a_{4}<i_{2}<i_{3}$ & 231 & 321 & 213 \\
				$ i_{1}<a_{4}<i_{2}<a_{5}<i_{3}$ & 132 & 231 & 132 \\
				$ i_{1}<a_{5}<i_{2}<a_{4}<i_{3}$ & 231 & 321 & 312 \\
				$ i_{1}<a_{4}<i_{2}<i_{3}<a_{5}$ & 123 & 231 & 132 \\
				$ i_{1}<a_{5}<i_{2}<i_{3}<a_{4}$ & 321 & 321 & 312 \\
				$\underline{ i_{1}<i_{2}<a_{4}<a_{5}<i_{3}}$ & 
				132 & 231 & 231 \\
				$\underline{ i_{1}<i_{2}<a_{5}<a_{4}<i_{3}}$ & 231 & 321 & 321 \\
				$\underline{ i_{1}<i_{2}<a_{4}<i_{3}<a_{5}}$ & 123 & 231 & 231 \\
				$\underline{ i_{1}<i_{2}<a_{5}<i_{3}<a_{4}}$ & 321 & 321 & 
				321 \\[0.07cm]
				\hline
				$a_{4}<a_{5}<i_{1}<i_{2}<i_{3}$ & 132 & 123 & 123 \\
				$a_{5}<a_{4}<i_{1}<i_{2}<i_{3}$ & 231 & 213 & 213 \\
				$a_{4}<i_{1}<i_{2}<i_{3}<a_{5}$ & 123 & 132 & 132 \\
				$a_{5}<i_{1}<i_{2}<i_{3}<a_{4}$ & 321 & 312 & 312 \\
				$\underline{i_{1}<i_{2}<i_{3}<a_{4}<a_{5}}$ & 
				\underline{213} & 231 & 231 \\
				$\underline{i_{1}<i_{2}<i_{3}<a_{5}<a_{4}}$ & 312 & 321 & 
				321 \\[0.07cm]
				\hline
			\end{tabular} 
			\caption{$M$-sets of cardinality $H_{52}(\lambda)
				|C_{\lambda}|/\fall{n}{5}$. \newline
				Underlines refer to Example~\ref{ex1}.}
			\label{mseth32}
		\end{table}
	}
	
}

\clearpage
\addcontentsline{toc}{section}{References}
\markright{References}{}

\end{document}